\documentclass[11pt,twoside]{article}
\usepackage{geometry}                		
\usepackage{graphicx}									
\usepackage{amssymb}\vspace{0.25cm}
\usepackage{mathtools}
\usepackage{setspace}
\usepackage{tikz-cd}
\usepackage[mathscr]{euscript}
\newcommand{\abs}[1]{\left\vert#1\right\vert}
\usepackage[normalem]{ulem}
\usepackage{manfnt}
\usepackage{pgfplots}
\usepackage{pifont}
\usepackage[safe]{tipa}
\usepackage{marvosym}
\usepackage{scalerel,stackengine}
\usepackage{titlesec}
\usepackage{makeidx}
\usepackage{centernot}

\usepackage{xspace}
\usepackage[noauto]{chappg}
\usepackage[OT2,T1]{fontenc}
\usepackage{hyperref}
\usepackage{amsmath}
\usepackage[amsmath,thmmarks]{ntheorem}

\usepackage{tabularx}
\newcolumntype{Y}{>{\raggedright\arraybackslash}X}

\newcommand\mA{%
$A$\xspace
}

\newcommand\mC{%
$C$\xspace
}

\newcommand\mE{%
$E$\xspace
}
\newcommand\mF{%
$F$\xspace
}
\newcommand\mG{%
$G$\xspace
}
\newcommand\mH{%
$H$\xspace
}

\newcommand\mK{%
$K$\xspace
}
\newcommand\mL{%
$L$\xspace
}
\newcommand\mM{%
$M$\xspace
}
\newcommand\mN{%
$N$\xspace
}

\newcommand\mP{%
$P$\xspace
}
\newcommand\mQ{%
$Q$\xspace
}

\newcommand\mX{%
$X$\xspace
}

\newcommand\mZ{%
$Z$\xspace
}
\newcommand\C{%
\mathbb{C}
}

\newcommand\Z{%
\mathbb{Z}
}
\newcommand\mfa{%
\mathfrak{a}
}
\newcommand\mfg{%
\mathfrak{g}
}
\newcommand\mfn{%
\mathfrak{n}
}
\newcommand\mfp{%
\mathfrak{p}
}
\newcommand\mfS{%
\mathfrak{S}
}

\newcommand\bPSL{%
\textbf{PSL}\xspace
}

\newcommand\ver{%
\text{ver}\hspace{0.05 cm}
}

\newcommand\bc{%
\textbf{c}\xspace
}
\newcommand\be{%
\textbf{e}\xspace
}
\newcommand\bn{%
\textbf{n}\xspace
}
\newcommand\bp{%
\textbf{p}\xspace
}
\newcommand\bCon{%
\textbf{Con}\xspace
}
\newcommand\bSp{%
\textbf{Sp}\xspace
}
\newcommand\bCl{%
\textbf{Cl}\xspace
}
\newcommand\bH{%
\textbf{H}\xspace
}
\newcommand\bK{%
\textbf{K}\xspace
}
\newcommand\bO{%
\textbf{O}\xspace
}
\newcommand\bQ{%
\textbf{Q}\xspace
}

\newcommand\bInd{%
\textbf{Ind}\xspace
}

\newcommand\bFnc{%
\textbf{Fnc}\hspace{0.05cm}
}
\newcommand\btr{%
\textbf{tr}\hspace{0.05cm}
}
\newcommand\sC{%
\mathcal{C}
}

\newcommand\sP{%
\mathcal{P}
}

\newcommand\Ad{%
\text{Ad}
}

\newcommand\Img{%
\text{Im}\hspace{0.05cm}
}

\newcommand\rank{%
\text{rank}\hspace{0.05cm}
}

\newcommand\balpha{%
\boldsymbol\alpha
}
\newcommand\bGamma{%
\boldsymbol\Gamma
}
\newcommand\ra{%
\rightarrow
}

\newcommand\vol{%
\text{vol}%\hspace{0.05cm}
}

\newcommand\dis{%
\text{dis}\hspace{0.05cm}
}
\newcommand\tr{%
\text{tr}\hspace{0.05cm}
}

\newcommand\restr[2]{%
{#1}|{#2}
}

\newcommand\vsx{%
\vphantom{\int_\int^\int}
}
\newcommand\vsy{%
\vphantom{\int}
}
\newcommand\hsx{%
\hspace{0.05cm}
}
\newcommand\hsy{%
\hspace{0.02cm}
}

\newcommand\rsC{%
\reflectbox{$\sC$}
\hspace{0.05cm}
}

\usepackage{scalerel,stackengine}
\stackMath
\newcommand\reallywidehat[1]{%
\savestack{\tmpbox}{\stretchto{%
  \scaleto{%
    \scalerel*[\widthof{\ensuremath{#1}}]{\kern-.6pt\bigwedge\kern-.6pt}%
    {\rule[-\textheight/2]{1ex}{\textheight}}%WIDTH-LIMITED BIG WEDGE
  }{\textheight}% 
}{0.5ex}}%
\stackon[1pt]{#1}{\tmpbox}%
}
\DeclareFontFamily{U}{wncy}{}
\DeclareFontShape{U}{wncy}{m}{n}{<->wncyr10}{}
\DeclareSymbolFont{mcy}{U}{wncy}{m}{n}
\DeclareMathSymbol{\Sh}{\mathord}{mcy}{"58} 

\makeatletter
\DeclareRobustCommand\widecheck[1]{{\mathpalette\@widecheck{#1}}}
\def\@widecheck#1#2{%
    \setbox\z@\hbox{\m@th$#1#2$}%
    \setbox\tw@\hbox{\m@th$#1%
       \widehat{%
          \vrule\@width\z@\@height\ht\z@
          \vrule\@height\z@\@width\wd\z@}$}%
    \dp\tw@-\ht\z@
    \@tempdima\ht\z@ \advance\@tempdima2\ht\tw@ \divide\@tempdima\thr@@
    \setbox\tw@\hbox{%
       \raise\@tempdima\hbox{\scalebox{1}[-1]{\lower\@tempdima\box
\tw@}}}%
    {\ooalign{\box\tw@ \cr \box\z@}}}
\makeatother

\newtheoremstyle{xx}% name of the style to be used
  {4pt}% measure of space to leave above the theorem. E.g.: 3pt
  {0pt}% measure of space to leave below the theorem. E.g.: 3pt
  {\upshape}% name of font to use in the body of the theorem
  {\bfseries}% name of head font
  {}% punctuation between head and body
  { }% space after theorem head; " " = normal interword space \footnotesize
  {}
  
\makeatletter 
 \newtheoremstyle{myu}%
  {\upshape\item[ \indent\indent\bf\underline{\theorem@headerfont ##2:}]}%
\makeatother
\makeatletter 
 \newtheoremstyle{myn}%
  {\item[\hskip\labelsep \ \bf ##1 \theorem@headerfont ##2.]}%
\makeatother
\theoremstyle{myn}
\theoremstyle{myu}
{\upshape}
\titleformat{\chapter}[display]
{\normalfont\filcenter\huge\bfseries}{}{0pt}{\large}
\usepackage[OT2,OT1]{fontenc} 
\newcommand\cyr
{
\renewcommand\rmdefault{wncyr} 
\renewcommand\sfdefault{wncyss} 
\renewcommand\encodingdefault{OT2} 
\normalfont
\selectfont
}

\DeclareTextFontCommand{\textcyr}{\cyr}

\makeindex 
\begin{document}

\title{\textbf{The Selberg Trace Formula IX:\\
Contribution from the Conjugacy Classes \\(The Regular Case)}}
\author{by\\
M. Scott Osborne$^*$\\
and\\
Garth Warner\footnote{ \ Research supported in part by the National Science Foundation}}

%\affil{University of Washington\\
%Seattle, Washington 98195}

\date{University of Washington\\
Seattle, Washington 98195}

\maketitle                              % Print title page.

\titlespacing*{\chapter}{0pt}{-50pt}{40pt}
\setlength{\parskip}{0.1em}
\include{_tocX}
\include{_Preface}
%\pagenumbering{bysection}
\setcounter{page}{1}

\setcounter{section}{1}
\renewcommand{\thepage}{1-\arabic{page}}
\setcounter{page}{1}
\setcounter{section}{-1}

\section
{
\qquad\qquad\qquad\qquad $\boldsymbol{\S}$\textbf{1}.\quad  Introduction
}
\setlength\parindent{2em}
\setcounter{theoremn}{0}

%%----------------------------------------------------------------------------------------------01
\ \indent 
This is the ninth in a projected series of papers in which we plan to come to grips with the Selberg trace formula, 
the ultimate objective being a reasonably explicit expression.  In our last publication [2-(h)], we isolated the contribution 
to the trace arising from the continuous spectrum, call it
\[
\bCon - \bSp(\alpha : \Gamma).
\]
Here, we shall initiate the study of the contribution to the trace arising from the conjugacy classes, call it 
\[
\bCon - \bCl(\alpha : \Gamma).
\]

Thus, in the usual notation, one has
\[
\tr(L_{G/\Gamma}^\dis(\alpha))\hsx = \hsx 
\bCon - \bCl(\alpha : \Gamma) + \bCon - \bSp(\alpha : \Gamma).
\]
And, thanks to Theorem 4.3 of [2-(h)], there is a formula for 
\[
\bCon - \bSp(\alpha : \Gamma)
\]
involving familiar ingredients, namely \bc-functions, exponentials, and \bInd-functions.  As for 
\[
\bCon - \bCl(\alpha : \Gamma),
\]
it will turn out that
\[
\bCon - \bCl(\alpha : \Gamma) \hsx = \hsx 
\sum\limits_{\sC_0}\sum\limits_{\sC \succeq \sC_0} \bCon - \bCl(\alpha : \Gamma_{\sC,\sC_0}),
\]
reflecting the partition
\[
\Gamma \hsx = \hsx  \coprod\limits_{\sC_0} \coprod\limits_{\sC \succeq \sC_0} \Gamma_{\sC,\sC_0}
\]
explained in \S 2.  The objective then is to find a formula for 
\[
\bCon - \bCl(\alpha : \Gamma_{\sC,\sC_0})
\]
involving familiar ingredients, namely orbital integrals (weighted or not) and their variants.

[Note: Strictly speaking (and in complete analogy with the conclusions of [2-(h)]), the formula will also contain a parameter 
\bH from the truncation space that, however, we shall ignore for the purpose of this Introduction.]

\vspace{0.1cm}

%%-----------------------------------------------------------------------------------------------------02

The simplest case is when $\sC = \sC_0$, $\Gamma_{\sC_0,\sC_0}$, being what we like to think of as the
$\sC_0$-regular elements of $\Gamma$.  
Each such element is necessarily semisimple and the identity component of its centralizer is contained in the 
$\Gamma$-Levi subgroup minimal with respect to ``$\in$'', although this need not be true of the full centralizer, 
a complicating circumstance.  
The situation is therefore similar to that considered by Arthur [1-(b), \S8], although a little more general due to the last mentioned point.  
The basis for our analysis will be the machinery set down in [2-(f)], applied to the case at hand.  
Its consistent use serves to smooth out most of the technical wrinkles, permitting the exposition to proceed in a systematic fashion.  

We begin in \S2 with a classification of the elements of $\Gamma$, introducing in particular the $\Gamma_{\sC,\sC_0}$.  
This material generalizes the rank-1 considerations of [2-(a), \S5].  
In \S3, the fine structure of $\Gamma_{\sC_0,\sC_0}$ is examined.   
\S4 is a brief exposition of the ``big picture'' and may be regarded as a supplement to this Introduction.  
The study of 
\[
\bCon - \bCl(\alpha : \Gamma_{\sC_0,\sC_0})
\]
is taken up in \S5, the main result being Theorem 5.4, which then, in \S6, is recast inductively as Theorem 6.3.

Finally, in what follows, the abbreviation TES will refer to our monograph, 
\textit{The Theory of Eisenstein Systems}, Academic Press, N.Y., 1981.

%%%%%%%%%%%%%%%%%%%%%%%%%%%%%%%%%%%%%%
%%%%%%%%%%%%%%%%%%%%%%%%%%%%%%%%%%%%%%
%%%%%%%%%%%%%%%%%%%%%%%%%%%%%%%%%%%%%%

\renewcommand{\thepage}{2-\arabic{page}}
\setcounter{page}{1}
\setcounter{section}{-1}

\section
{
\qquad\qquad $\boldsymbol{\S}$\textbf{2}.\quad  Classification of the Elements of $\Gamma$
}
\setlength\parindent{2em}
\setcounter{theoremn}{0}
%%----------------------------------------------------------------------------------------------01
\ \indent 
The purpose of this \S \  is to devise a delineation of the elements of $\Gamma$ suitable for the calculations which are to follow 
(here and elsewhere).

Given $\gamma \in \Gamma$, let us agree to write $\{\gamma\}_G$ 
(respectively $\{\gamma\}_\Gamma$) for its $G$-conjugacy class
(respectively $\Gamma$-conjugacy class), $G_\gamma$ (respectively $\Gamma_\gamma$) for its $G$-centralizer 
(respectively $\Gamma$-centralizer).

Put
\[
Z_\Gamma \hsx = \hsx \text{center of} \ \Gamma.
\]

\vspace{0.1cm}

{\small\bf Proposition 2.1} \ %01
Suppose that $\gamma \in Z_\Gamma$ $-$then $\gamma$ is semisimple.  
Furthermore, $\Gamma_\gamma$ is a nonuniform lattice in $G_\gamma$.

\vspace{0.1cm}

[To prove this, one need only repeat the discussion on p. 19 of [2-(a)], the assumption there that 
$\rank(\Gamma) = 1$ being of no relevance at all.]

\vspace{0.3cm}

An element $\gamma \in \Gamma$ is said to be $G$-regular provided that $\gamma$ belongs to no proper \linebreak
$\Gamma$-cuspidal parabolic subgroup of \mG.  
Denote by $\Gamma_G$ the set of such $-$then $\Gamma_G$ is obviously invariant under $\Gamma$-conjugacy.  
Calling $[\Gamma_G]$ some choice of representatives for the $\Gamma$-conjugacy classes in $\Gamma_G$, form 
\[
\coprod\limits_{\gamma \in [\Gamma_G]} \ (G / \Gamma_\gamma) \times \{\gamma\}
\]
and let $\Phi$ be the canonical map from this set to \mG, viz.
\[
\Phi(x \Gamma_\gamma, \gamma) \hsx = \hsx x \gamma x^{-1}.
\]

\vspace{0.1cm}

{\small\bf Proposition 2.2} \ %02
$\Phi$ is a proper map.

\vspace{0.1cm}

We shall need  a preliminary result.

\vspace{0.2cm}

{\small\bf Lemma 2.3} \ %03
Let \mP be a $\Gamma$-percuspidal parabolic subgroup of \mG with split component \mA.  Suppose that
\[
\begin{cases}
\ \{a_n\} \in A[t]\\
\ \{\gamma_n\} \in \Gamma
\end{cases}
\]
are sequences such that
\[
%a_n^\lambda \ra -\infty \qquad (\lambda \in \sum_P^0 \ (\mfg,\mfa))
a_n^\lambda \ra -\infty \qquad (\lambda \in \raisebox{-0.075cm}{\scalebox{1.75}{$\Sigma$}}_P^0 (\mfg,\mfa))
\]
%%----------------------------------------------------------------------------------------------02
and 
\[
a_n \gamma_n a_n^{-1} \qquad \text{stays bounded}.
\]
Then eventually
\[
\gamma_n \in P_F, \quad F = \raisebox{-0.1cm}{\scalebox{1.75}{$\Sigma$}}_P^0  (\mfg,\mfa) - \{\lambda\}.
\]

\vspace{0.1cm}

{\small\bf Proof.}
Let $\eta \in \Gamma \cap N_F$ $-$then 
\[
a_n \hsy \eta \hsy a_n^{-1} \ra 1
\]
and so
\[
(a_n \hsy \gamma_n \hsy a_n^{-1}) \hsy
a_n \hsy \eta \hsy  a_n^{-1}
(a_n \hsy \gamma_n\hsy  a_n^{-1})^{-1}
\ra 1
\]
or still
\[
a_n \hsy \gamma_n \hsy \eta \hsy \gamma_n^{-1} a_n^{-1} \ra 1.
\]
Thanks to the lemma on p. 47 of TES, 
\[
\gamma_n \hsy\eta \hsy\gamma_n^{-1} \in N
\]
eventually.   
Consequently, upon taking generators for $\Gamma \cap N_F$, 
\[
\gamma_n(\Gamma \cap N_F) \gamma_n^{-1} \subset N
\]
eventually.  This implies that for all sufficiently large $n$, 
\[
\Ad(\gamma_n) \mfn_F \subset \mfn,
\]
from which, passing to orthocomplements, 
\[
\Ad(\gamma_n) \mfp_F \supset \mfp,
\]
i.e., 
\[
\gamma_n P_F \gamma_n^{-1} \supset P.
\]
Therefore $\forall \ n \gg 0$, 
\[
\gamma_n P_F \gamma_n^{-1} \hsx = \hsx P_F \implies \gamma_n \in P_F,
\]
as desired.

\vspace{0.3cm}

{\small\bf Proof of Proposition 2.2}
Let $C_G$ be a compact subset of \mG $-$then we must show that $\Phi^{-1}(C_G)$ is a compact subset of
\[
\coprod\limits_{\gamma \in [\Gamma_G]} \ (G / \Gamma_\gamma) \times \{\gamma\}.
\]
%%----------------------------------------------------------------------------------------------03
For this purpose, let 
\[
\{(x_n \Gamma_{\gamma_n}, \gamma_n)\}
\]
be a sequence in $\Phi(C_G)^{-1}$.  
Since \mG can be covered by finitely many $\mfS \bullet \Gamma$, 
$\mfS$ a Siegel domain relative to a $\Gamma$-percuspidal \mP, 
there is no loss of generality in supposing to begin with that $x_n \in \mfS \bullet \Gamma \  \forall \ n$.  
That being, write as usual 
\[
x_n \hsx = \hsx k_n \hsy a_n \hsy s_n \delta_n.
\]
Because 
\[
k_n \hsy a_n \hsy s_n \hsy a_n^{-1}
\]
stays bounded, the same holds for 
\[
a_n \delta_n \gamma_n \delta_n^{-1} a_n^{-1}.
\]
But 
\[
\delta_n \gamma_n \delta_n^{-1} \notin P_F \qquad \forall \ F \neq \Sigma_P^0(\mfg,\mfa)
\]
and so, on the basis of the foregoing lemma, $a_n$ must stay bounded too.  
Therefore, after passing to a subsequence if necessary, we may assume that $k_n$, $a_n$, and $s_n$ are all convergent, 
along with $x_n \gamma_n x_n^{-1}$.  
Accordingly, $\delta_n \gamma_n \delta_n^{-1}$ is convergent, hence is eventually constant.  
Thus, by definition of $[\Gamma_G]$, $\forall \ n \gg 0$, 
$\gamma_n = \gamma$ and $\delta_n \Gamma_\gamma = \delta \Gamma$, implying that $\lim x_n \Gamma_\gamma$ exists.  
This establishes the compactness in 
\[
\coprod\limits_{\gamma \in [\Gamma_G]} \ (G / \Gamma_\gamma) \times \{\gamma\}
\]
of $\Phi^{-1}(C_G)$.

\vspace{0.3cm}

The following points are immediate corollaries.

(1) 
\quad $\forall \ \gamma \in \Gamma_G$, $\gamma$ is semisimple and $\Gamma_\gamma$ is a uniform lattice in $G_\gamma$.

(2)
\quad $\forall$ compactum $C_G \subset G$, 
\[
\#\big(\{\gamma \in [\Gamma_G] : \{\gamma\}_G \cap C_G \neq \emptyset\}\big) \hsx < \hsx +\infty.
\]

(3)
\quad $\forall$ compactum $C_G \subset G$ and $\forall \ \gamma \in \Gamma_G$,
\[
\{x \in G : x \gamma x^{-1} \in C_G\}
\]
is compact mod $G_\gamma$.
%%----------------------------------------------------------------------------------------------04

(4) 
\quad$\forall \ \gamma \in \Gamma_G$, 
\[
\#\bigg(\{\gamma^\prime \in [\Gamma_G] : \{\gamma^\prime\}_G = \{\gamma\}_G \}\bigg) \hsx < \hsx +\infty.
\]

Observe that no $\Gamma$-central element can belong to $\Gamma_G$.  
Indeed, if $\gamma \in Z_\Gamma$, then $\Gamma_\gamma$ is a nonuniform lattice in $G_\gamma$, 
whereas if $\gamma \in \Gamma_G$, then $\Gamma_\gamma$ is a uniform lattice in $G_\gamma$.

Suppose now that \mP is a $\Gamma$-cuspidal parabolic subgroup of \mG $-$then an element $\gamma \in \Gamma$ 
is said to be \mP-regular if $\gamma \in P$ but $\gamma \notin P^\prime$ for all $P^\prime \prec P$.  
Denote by $\Gamma_P$ the set of such $-$then
\[
\Gamma \hsx = \hsx \bigcup \Gamma_P,
\]
although, of course, there is overlap in the union on the right.

\vspace{0.2cm}

{\small\bf Proposition 2.4} \ %04
If 
\[
\Gamma_{P_1} \cap \Gamma_{P_2} \hsx \neq \hsx \emptyset, 
\]
then $P_1$ and $P_2$ are associate.

\vspace{0.1cm}

Admitting this momentarily, given an association class $\sC$, put
\[
\Gamma_\sC \hsx = \hsx \bigcup\limits_{P \in \sC} \Gamma_P.
\]
Since 
\[
\gamma \Gamma_P \gamma^{-1} \hsx = \hsx \Gamma_{\gamma P \gamma^{-1}},
\]
it is clear that $\Gamma_\sC$ is invariant under $\Gamma$-conjugacy and, by the above, 
\[
\sC^\prime \hsx \neq \hsx \sC^{\prime\prime}
\implies 
\Gamma_{\sC^\prime} \cap \Gamma_{\sC^{\prime\prime}} \hsx = \hsx \emptyset,
\]
so
\[
\Gamma \hsx = \hsx \coprod\limits_\sC \Gamma_\sC.
\]
Needless to say, 
\[
\Gamma_{\{G\}} \hsx = \hsx \Gamma_G.
\]

In general, if \mP is a $\Gamma$-cuspidal parabolic subgroup of \mG with unipotent radical, \mN, 
then a Levi subgroup \mL of \mP is a closed reductive subgroup with the property that the multiplication 
$L \times N \ra P$ is an isomrphism of analytic manifolds (cf. TES, p. 31).  
To reflect $\Gamma$-cuspidality, it will be best to specialize this notion, using the
\[
\begin{cases}
\ (G^*,\Gamma^*) \\ 
\ (G^\#,\Gamma^\#) 
\end{cases}
\]
%%----------------------------------------------------------------------------------------------05
formalism in TES (pp. 40-41), putting $1 \leq i \leq r$
\[
\begin{cases}
\ G_i^\# \ \hsx = \hsx G^\# / \prod\limits_{j \neq i} G_j^*\\
\ \Gamma_i^\#\ \hsx = \hsx \Gamma^\# \bullet \prod\limits_{j \neq i} G_j^* / \prod\limits_{j \neq i} G_j^*
\end{cases}
\]
and 
\[
x_i^\# \hsx = \hsx (x Z \bullet G_{c,ss}) \bullet \prod\limits_{j \neq i} G_j^* \qquad (x \in G).
\]
\vspace{0.1cm}
By a $\Gamma$-Levi subgroup \mL of \mP we shall then understand a Levi subgroup \mL of \mP such that 

(i) \quad When $i \leq r_1$, $L_i^\# = G_i^\#$.

(ii) \quad When $r_1 < i \leq r_2$, $L_i^\#$ is a Levi subgroup of $P_i^\#$ (per $\Gamma_i^\#$).

(iii) \quad When $r_2 < i \leq r$, $\Ad(L_i^\#)$ is a \bQ-Levi subgroup of $\Ad(P_i^\#)$ per  $\Ad(\Gamma_i^\#)$).

Generically, let 
\[
\gamma \hsx = \hsx  \gamma_s \gamma_u
\]
be the Jordan decomposition of $\gamma$.  
Note that 
\[
\gamma_s \in P \quad \text{iff} \quad \forall \ i,\  (\gamma_s^\#)_i \in P_i^\#.
\]
Moreover, if \mL is a $\Gamma$-Levi subgroup of \mP, then 
\[
\gamma_s \in L \quad \text{iff} \quad \forall \ i, \  (\gamma_s^\#)_i \in L_i^\#. 
\]

\vspace{0.3cm}

{\small\bf Lemma 2.5} \ %05
Suppose that $\gamma \in \Gamma_P$ $-$then there exists a $\Gamma$-Levi subgroup $L = M \bullet A$ of \mP such that 
$\gamma_s \in M$.  
In addition, \mP is minimal with respect to ``$\gamma_s \in P$''.

\vspace{0.1cm}

{\small\bf {Proof.}} \ 
In view of the reductions outlined above, we can deal with each possibility separately, mentally making changes in the notation.  
Start off with any \mM and let $\{\delta\} = M \cap \gamma N$ $-$then $\delta \in \Gamma_M$ and we claim that $\delta$ is 
\mM-regular, hence semisimple.  
For otherwise, 
\[
\delta \in \prescript{\prime}{}P \implies \delta N \subset P^\prime \implies \gamma \in P^\prime,  
\]
an impossibility.  Let 
\[
\mfn(\delta) \hsx = \hsx \Img(\Ad(\delta) - 1)
\]
and use the surjection
\[
\phi : 
\begin{cases}
\  \mfn(\delta) \times N_\delta \ra N \\
\ (X,n) \ra \delta^{-1} \exp(X) (\delta n) \exp(-X)
\end{cases}
\]
%%----------------------------------------------------------------------------------------------06
to write
\[
\gamma \hsx = \hsx \delta \phi(X,n) \hsx = \hsx \exp(X) (\delta n) \exp(-X).
\]
Since $\delta$ and $n$ commute, 
\[
\gamma \hsx = \hsx [\exp(X)  \delta \exp(-X)] \bullet  [\exp(X)  n \exp(-X)] 
\]
is the Jordan decomposition of $\gamma$.  In particular,
\[
\gamma_s \in \exp(X)  M \exp(-X),
\]
leading to the first assertion.  
As for the second, if the contrary were true, then the \mM-regularity of $\gamma_s$ would be violated.

\vspace{0.3cm}

In passing, observe that \mL is actually minimal with respect to ``$\gamma_s \in L$''.

\vspace{0.3cm}

{\small\bf {Proof of Proposition 2.4.}} \ 
Let 
\[
\gamma \in \Gamma_{P_1} \cap \Gamma_{P_2}.
\]
To prove that $P_1$ and $P_2$ are associate, we need only produce an 
$L = M \bullet A$ with $\gamma_s \in M$ that is simultaneously a $\Gamma$-Levi subgroup for both $P_1$ and $P_2$.  
To this end, we shall proceed on a case-by-case basis.  
If $i \leq r_1$, then there is nothing to prove.  
If $r_1 < i \leq r_2$, then either $P_1 = P_2$ or $P_1 \neq P_2$ and in the latter situation
\[
P_1 \cap P_2 
\hsx = \hsx 
L 
\hsx = \hsx  
M \bullet A 
\implies 
\gamma \in M 
\implies 
\gamma 
\hsx = \hsx 
\gamma_s.
\]
Finally, if $r_2 < i \leq r$, then \mL is minimal with respect to ``$\gamma_s \in L$'' iff \mA is maximal with respect to 
$A \subset G_{\gamma_s}^0$, meaning that \mA is a maximal \bQ-split torus in $G_{\gamma_s}^0$.  
Any two such are conjugate.  
Because 
\[
P_1 \cap G_{\gamma_s}^0 \quad \text{and} \quad P_2 \cap G_{\gamma_s}^0
\]
are minimal \bQ-parabolic subgroups, they share a split component and we can take for \mL its \mG-centralizer.

\vspace{0.3cm}

Having decomposed $\Gamma$ as a disjoint union 
\[
\coprod\limits_{\sC} \ \Gamma_\sC,
\]
the next step will be to fix a $\sC_0$ and decompose $\Gamma_{\sC_0}$ still further in terms of 
$\sC \succeq \sC_0$.

\vspace{0.3cm}

%%----------------------------------------------------------------------------------------------07
{\small\bf Lemma 2.6} \ %07
Let $L = M \bullet A$ be a $\Gamma$-Levi subgroup of \mP.  
Suppose that
\[
\gamma \in P \quad \text{and} \quad \gamma_s \in M.
\]
Then the following conditions are equivalent:
%\[
\begin{align*}
\ \text{(i)} \ N_\gamma \hsx = \hsx \{1\}; \qquad \hspace{0.25cm}  &\text{(iii)} \ P_{\gamma_s} \subset L; \hspace{0.25cm} \\
\ \text{(ii)} \ N_{\gamma_s} \hsx = \hsx \{1\}; \qquad  \ &\text{(iv)} \ G_{\gamma_s}^0 \subset L.
\end{align*}
%\]

[The verification is straightforward, hence can be omitted.]

\vspace{0.3cm}

{\small\bf Remark.} \ %0
As can be seen by example, one cannot improve (iv) to read $G_{\gamma_s} \subset L$.  
For instance, take
\[
G \hsx = \hsx \bPSL_2(\C), \quad \Gamma  \hsx = \hsx \bPSL_2\bigl(\Z[\sqrt{-1}]\bigr)
\]
and consider
\[
\gamma \ = \ 
\begin{pmatrix}
0 &\pm \sqrt{-1}\\
\mp \sqrt{-1} &0\\
\end{pmatrix}
.
\]

\vspace{0.2cm}

{\small\bf Proposition 2.7} \ %08
If $P_1$ and $P_2$ are minimal with respect to ``$\gamma \in P$ and $N_\gamma = \{1\}$'', then 
$P_1$ and $P_2$ are associate.

\vspace{0.1cm}

{\small\bf Proof.} \ 
Once again, we shall examine cases.  
If $i \leq r_1$, then there is nothing to prove.  
If $r_1 < i \leq r_2$, then either $\gamma$ and $G_{\gamma_s}^0$ are contained in a proper $\Gamma$-cuspidal or they are not.  
Finally, if $r_2 < i \leq r$, then $\gamma_s$ lies in $L_\bQ$ and 
\[
\begin{cases}
\ G_{\gamma_s}^0 \subset L_1 \\
\ G_{\gamma_s}^0 \subset L_2
\end{cases}
.
\]
Denote by \mH the reductive algebraic \bQ-group generated by $\gamma_s$ and $G_{\gamma_s}^0$ $-$then 
$H \subset G_{\gamma_s}$ and $L_1$ and $L_2$ are minimal with respect to ``$H \subset L$'', thus are the 
centralizers of maximal \bQ-split tori in $C_G(H)^0$ (which is reductive), so are conjugate.

\vspace{0.3cm}

Fix an association class $\sC_0$.  Given a $\sC \succeq \sC_0$, call
\[
\Gamma_{\sC,\sC_0}
\]
the set of $\gamma \in \Gamma_{\sC_0}$ for which there exists a $P \in \sC$ minimal with respect to 
``$\gamma \in P$ and $N_\gamma = \{1\}$''.  
On the basis of the proposition supra,
\[
\sC^\prime \hsx \neq \hsx \sC^{\prime\prime} 
\implies 
\Gamma_{\sC^\prime,\sC_0} \cap \Gamma_{\sC^{\prime\prime}, \sC_0} \hsx = \hsx \emptyset,
\]
%%----------------------------------------------------------------------------------------------08
hence 
\[
\Gamma_{\sC_0} \hsx = \hsx \coprod \limits_{\sC \succeq \sC_0} \Gamma_{\sC,\sC_0}.
\]
In particular, 
\[
Z_\Gamma \subset \Gamma_{\{G\},\sC_0},
\]
$\sC_0$ the association class made up of the $\Gamma$-percuspidals.

Heuristically, 
\[
\Gamma_{\sC_0,\sC_0}
\]
are the $\sC_0$-regular elements of $\Gamma$ (cf. [2-(a)] when $\rank (\Gamma) = 1$).

\vspace{0.3cm}

{\small\bf Remark.} \ %0
Because 
\[
\begin{cases}
\ \gamma \Gamma_{\sC,\sC_0} \gamma^{-1} = \Gamma_{\sC,\sC_0} \qquad (\gamma \in \Gamma)\\
\ \gamma \in \Gamma_{\sC,\sC_0} \cap P \implies \gamma \bullet (\Gamma \cap N) \subset \Gamma_{\sC,\sC_0}
\end{cases}
,
\]
the partition 
\[
\Gamma \hsx = \hsx \coprod \limits_{\sC_0} \coprod \limits_{\sC \succeq \sC_0} \Gamma_{\sC,\sC_0} 
\]
satisfies the general assumptions laid down on pp. 1376-1377 of [2-(f)].  
The associated $\Gamma$-compatible families are therefore estimable (cf. [2-(f), p. 1413]).

\vspace{0.1cm}
%%%%%%%%%%%%%%%%%%%%%%%%%%%%%%%%%%%%%%
%%%%%%%%%%%%%%%%%%%%%%%%%%%%%%%%%%%%%%
%%%%%%%%%%%%%%%%%%%%%%%%%%%%%%%%%%%%%%

\renewcommand{\thepage}{3-\arabic{page}}
\setcounter{page}{1}
\setcounter{section}{-1}

\section
{
\qquad\qquad\qquad\qquad $\boldsymbol{\S}$\textbf{3}.\quad  The Structure of $\Gamma_{\sC_0,\sC_0}$
}
\setlength\parindent{2em}
\setcounter{theoremn}{0}
%%----------------------------------------------------------------------------------------------01
\ \indent 

The purpose of this \S\  is to look more closely at the structure of $\Gamma_{\sC_0,\sC_0}$, the results obtained being essential preparation for the analysis carried out in \S 5 infra.

We shall begin with a series of simple observations.

Put $\Delta = \Gamma_{\sC_0,\sC_0}$ and decompose $\Delta$ into $\Gamma$-conjugacy classes $\Delta_i$: 
\[
\Delta \ = \ \coprod \Delta_i.
\]
Fix an index $i$, say $i = 0$, and fix a $\gamma_0 \in \Delta_0$.

Let $P_0$ be an element of $\sC_0$ minimal with respect to ``$\gamma_0 \in P$ and $N_{\gamma_0} = \{1\}$''.  
Choose, as is possible (cf. Lemma 2.5), a $\Gamma$-Levi subgroup $L_0$ of $P_0$ such that 
$ss(\gamma_0) \in L_0$ $-$then the centralizer of $ss(\gamma_0)$ in $N_0$ is necessarily trivial (cf. Lemma 2.6), 
thus actually $\gamma_0 = ss(\gamma_0)$ and so $\gamma_0$ is semisimple.

Accordingly, $L_0$ is minimal with respect to ``$\gamma_0 \in L$'' and is in fact the only $\Gamma$-Levi subgroup with this property.  
To see this, examine cases, the nontrivial one being the algebraic situation.  
But, from the proof of Proposition 2.4, $L_0$ is minimal with respect to ``$\gamma_0 \in L$'' iff $L_0$ is the 
$G$-centralizer of a maximal \bQ-split torus in $G_{\gamma_0}^0$, itself unique as $G_{\gamma_0}^0 \subset L_0$ 
(cf. Lemma 2.6), from which the uniqueness of $L_0$.  
Of course, there is exactly one splitting $L_0 = M_0 \bullet A_0$ with $\gamma_0 \in M_0$, 

Denote by 
\[
\sP_0(L_0)
\]
the set of parabolics having $L_0$ for a Levi factor $-$then all such are $\Gamma$-cuspidal.  
To see this, again examine cases, the nontrivial one being the rank-1 situation, covered fortunately, by Lemma 9.3 in [2-(a)].  
Obviously, therefore, 
\[
\sP_0(L_0) \ = \ \{P_0 \in \sC_0 : \gamma_0 \in P_0\}.
\]

If now by 
\[
\sP_0(L_0)
\]
we understand
\[
\{P \in \sC_\Gamma : \gamma_0 \in P\},
\]
then it is clear that
\[
\sP_0(L_0)  \ = \ \big\{P \in \sC_\Gamma : P \succeq P_0 \ \bigl(\exists \ P_0 \in \sP_0(L_0) \bigr) \big\}.
\]
%%----------------------------------------------------------------------------------------------02
Evidently, 
\[
\forall \ P \in \sP(L_0) : \Gamma_{\gamma_0} \cap P \hsx = \hsx \Gamma_{\gamma_0} \cap L,
\]
$L \supset L_0$ is a $\Gamma$-Levi subgroup of \mP.

Given
\[
\begin{cases}
\ P \in \sP(L_0)\\
Q \in \sC_\Gamma
\end{cases}
,
\]
let
\[
\Delta_0(P)_Q \ = \ \{\gamma \gamma_0 \gamma^{-1} : \gamma P \gamma^{-1} \ = \ Q\}.
\]
Then
\[
\Delta_0(P)_Q \hsx \neq \hsx \emptyset
\]
iff \mP and \mQ are $\Gamma$-conjugate.  
Moreover, 
\[
\Delta_0(P)_P
\]
is the $\Gamma \cap P$-conjugacy class of $\gamma_0$.

\vspace{0.3cm}

{\small\bf Lemma 3.1} \ %3.1
$\forall \ \gamma \in \Gamma$, 
\[
\Delta_0(P)_{\gamma Q \gamma^{-1}} \ = \ \gamma \Delta_0(P)_Q \gamma^{-1}.
\]

\qquad [The verification is immediate.]

Let \mQ be an element of $\sC_\Gamma$ with the property that 
\[
\Delta_0 \cap Q \hsx \neq \hsx \emptyset.
\]
Let \mE be the set of all $P \in \sP(L_0)$ that are $\Gamma$-conjugate to \mQ $-$then we claim that \mE is a 
$\Gamma$-conjugacy class in $\sP(L_0)$.  
To check that this is the case, we need only convince ourselves that a $\Gamma$-conjugate of \mQ is in $\sP(L_0)$.   
But, for some $\gamma \in \Gamma$, 
\[
\gamma \gamma_0 \gamma^{-1} \in \Delta_0 \cap Q
\]
\[
\implies \gamma_0 \in \gamma^{-1} Q \gamma \implies \gamma^{-1} Q \gamma \in \sP(L_0),
\]
as desired.  Consequently, 
\[
\Delta_0 \cap Q \hsx \subset \hsx \bigcup\limits_{P \in E} \Delta_0(P)_Q.
\]
On the other hand, from the definitions,
\[
\forall \ P \in E, \quad \Delta_0(P)_Q \hsx \subset \hsx \Delta_0 \cap Q.
\]
%%----------------------------------------------------------------------------------------------03
Therefore
\[
\Delta_0 \cap Q \hsx = \hsx \bigcup\limits_{P \in E} \Delta_0(P)_Q.
\]
To render this union disjoint, let $E(\gamma_0)$ stand for a set of representatives from \mE per 
$\Gamma_{\gamma_0}$-conjugacy $-$then the lemma infra implies that
\[
\Delta_0 \cap Q \hsx = \hsx \coprod\limits_{P \in E(\gamma_0)} \Delta_0(P)_Q.
\]

\vspace{0.3cm}

{\small\bf Lemma 3.2} \ %3.2
Suppose that 
\[
P_1, P_2 \in \sP(L_0)
\]
are $\Gamma$-conjugate to a $Q \in \sC_\Gamma$ $-$then the following are equivalent:
\begin{align*}
&\vsy(i) \quad\ \  \Delta_0(P_1)_Q \hsx = \hsx \Delta_0(P_2)_Q;\\
&\vsy(ii) \quad\  \Delta_0(P_1)_Q \hsx \cap \hsx \Delta_0(P_2)_Q \neq \emptyset;\\
&\vsy(iii) \quad P_1 \ \text{and} \ P_2  \ \text{are $\Gamma_{\gamma_0}$-conjugate}.
\end{align*}

[The verification is immediate.]

\vspace{0.3cm}

Up to this point, we have worked with a fixed $\gamma_0$.  
For the remainder of this \S, it will be necessary to work with a variable $\gamma$, 
the corresponding notational changes being 
\[
M_0 \bullet A_0 \ = \ L_0 \ra L_0(\gamma) \ = \ M_0(\gamma) \bullet A_0(\gamma)
\]
\[
\begin{cases}
\ \sP_0(L_0) \ra \sP_0(L_0(\gamma)) \qquad P_0(\gamma)\\
\ \sP(L_0) \ \ra \sP(L_0(\gamma)) \qquad \  P(\gamma)
\end{cases}
\]
to signify the dependence on $\gamma$.

\vspace{0.3cm}

{\small\bf Proposition 3.3} \ %3.3
Let \mC be a compact subset of \mG $-$then 
\[
\#\big(\bigl\{ \{ \gamma\}_\Gamma : \gamma \in \Gamma_{\sC_0,\sC_0} 
\ \& \ 
\{\gamma\}_G \cap C \neq \emptyset \bigr\} \big ) \ < \ +\infty.
\]

[Note: \ 
When rank $(\Gamma) = 1$, this is Proposition 5.15 in [2-(a)].]

\vspace{0.1cm}

To get at this, some preparation will be required.  
We can certainly assume that $K \bullet C \bullet K = C$.  
Given now a $P \in \sC_\Gamma$, put 
\[
C_P \ = \ C \bullet N \cap M \qquad (\sim (C \cap S) \bullet N / N),
\]
%%----------------------------------------------------------------------------------------------04
a compact subset of \mM.

\vspace{0.3cm}

{\small\bf Lemma 3.4} \ %3.4
Let $\gamma \in \Gamma_{\sC_0,\sC_0}$.  
Suppose that $P_0(\gamma) \in \sP_0(L_0(\gamma))$ $-$then the following are equivalent:
\begin{align*}
&\vsy(i) \quad\  \{\gamma\}_G \cap C \neq \emptyset;\\
&\vsy(ii) \quad \{\gamma\}_{M_0(\gamma)} \cap C_{P_0(\gamma)} \neq \emptyset.
\end{align*}

[Note: To avoid any possibility of confusion, the relevant Langlands decomposition of $P_0(\gamma)$ is 
$M_0(\gamma) \bullet A_0(\gamma) \bullet N_0(\gamma)$, where $\gamma \in M_0(\gamma)$.]\\
{\small\bf Proof.} \ 
$(i) \ \implies \ (ii)$.  
By hypothesis, there exists an $x \in G$ such that $x \gamma x^{-1} \in C$.  
Per $P_0(\gamma)$, write, as usual, $x = k m n a$.  
Because $a \gamma a^{-1} = \gamma$, we have
\begin{align*}
&\vsx \quad x \gamma x^{-1} = \ k m n \gamma n^{-1} m^{-1} k^{-1} \in C\\
&\vsy \implies m n \gamma n^{-1} m^{-1} \in K \bullet C \bullet K = C\\
&\vsy \implies m \gamma m^{-1} \bullet m(\gamma^{-1} n \gamma n^{-1}) m^{-1}| \in C\\
&\vsy \implies m \gamma m^{-1} \in C \bullet N_0(\gamma)\\
&\vsy \implies \{\gamma\}_{M_0(\gamma)} \cap C_{P_0(\gamma)} \neq \emptyset.
\end{align*}

$(ii) \ \implies \ (i)$.  
By hypothesis, there exists an $m \in M_0(\gamma)$ such that $m \gamma m^{-1} \in C_{P_0(\gamma)}$.  
Accordingly, $m \gamma m^{-1} = cn$.  
Since
\[
\det (\Ad (m \gamma m^{-1} | \mfn_0 (\gamma) - 1) \neq 0, 
\]
there exists an $n_0 \in N_0(\gamma)$ such that 
\[
n_0(m \gamma m^{-1}) n_0^{-1} \ = \ m \gamma m^{-1} n^{-1} \in C
\]
\[
\implies \{\gamma\}_G \cap C \neq \emptyset.
\]

{\small\bf Remark.} \ 
The second criterion of the lemma does not depend on the choice of $P_0(\gamma) \in \sP_0(L_0(\gamma))$.  
Moreover, the use of $M_0(\gamma)$ is not crucial: One can move to the ambient special ``$M$'' provided $\gamma$ is replaced by $\delta$, $\{\delta\} = $ ``$M$'' $\cap \hsx\hsx \gamma N_0(\gamma)$.

\vspace{0.3cm}

Passing to the proof of the proposition, let us proceed by contradiction and assume that the cardinality in question is infinite.  
There will then be infinitely many $\{\gamma\}_\Gamma$ $(\gamma \in \Gamma_{\sC_0,\sC_0})$ having a member in some fixed $P_0(\gamma) \equiv P$ (for short), this because there are but finitely many $P_{i \mu}$.  
Working with the special ``$M$'', as always from $\gamma \in \Gamma$ we determine $\delta \in \Gamma_M$ via 
$\{\delta \} = M \cap \gamma N$, the $\Gamma$-conjugacy class of $\gamma$ filling out the $\Gamma_M$-conjugacy class 
%%----------------------------------------------------------------------------------------------05
of $\delta$.  And, in the case at hand, $\delta$ is $M$-regular.  
The number of $M$-conjugacy classes of $M$-regular elements that can meet $C_P$ is finite (cf. \S 2).  
So, there are infinitely many 
$\{\gamma\}_\Gamma$ $(\gamma \in \Gamma_{\sC_0,\sC_0})$
producing a fixed $\delta$, the corresponding set of $\gamma$ being precisely $\Gamma \cap \delta N$.  
However, we claim that $\Gamma \cap \delta N$ is the union of 
\[
\abs{\det (\Ad (\delta) | \mfn - 1)}
\]
$\Gamma \cap N$-conjugacy classes.  
Granted this, we have our contradiction.

The plan will be to appeal to a generality from [2-(a)], namely:

\vspace{0.1cm}

{\small\bf Proposition.} \ 
Let \mN be a connected, simply connected nilpotent Lie group; let $\Gamma$ be a lattice in \mN.  
Suppose that $\phi:N \ra N$ is an automorphism of \mN carrying $\Gamma$ into itself with 
$\det (d \phi - 1) \neq 0$.  
Put
\[
\Delta \ = \ \{n \in N : \phi(n) n^{-1} \in \Gamma\}.
\]
Then $\Delta$ is a (finite) union of $\abs{\det (d \phi - 1)}$ right cosets of $\Gamma$.

[For details, the reader is referred to pp. 24-25 of [2-(a)].  
The argument, as given there, is marginally incomplete, so we shall take this opportunity to set things straight.  
The difficulty is that in general the Leray-Serre spectral sequence uses local coefficients when the base space is not simply connected, a point that we had overlooked at the time.  
But here the local coefficients of $H_{\raisebox{-0.1cm}{\textbf{*}}}$ of the fiber are global.  
Thus, take a path in $Z \bullet \Gamma \backslash N$ and lift it to a path joining 1 and $z \gamma$ $-$then we must show that multiplication by $z \gamma$ on 
\[
(Z \cap \Gamma) \backslash Z \hsx \sim \hsx \Gamma \backslash Z \bullet \Gamma
\]
is the identity map on 
\[
H_{\raisebox{-0.1cm}{\textbf{*}}} ((Z \cap \Gamma) \backslash Z).
\]
Because \mZ is the center of \mN, multiplication by $\gamma$ is the identity and, as \mZ is arcwise connected, multiplication by $z$ is homotopic to the identity.]

We come now to the claim.  
Simply fix a $\gamma \in \Gamma \cap \delta N$ and, in the notation employed above, let
\[
\phi, n \ra \gamma^{-1} \hsy n \hsy \gamma.
\]
Conclude by observing that 
\begin{align*}
\det ( \Ad (\delta) | \mfn - 1) \ 
&= \ \det ( \Ad (\gamma) | \mfn - 1)\\
&= \det ( \Ad (\gamma) | \mfn) \bullet \det (1 - \Ad (\gamma^{-1} | \mfn), 
\end{align*}
with 
\[
\abs{\det ( \Ad (\gamma) | \mfn)} \ = \  1.
\]

%%%%%%%%%%%%%%%%%%%%%%%%%%%%%%%%%%%%%%
%%%%%%%%%%%%%%%%%%%%%%%%%%%%%%%%%%%%%%
%%%%%%%%%%%%%%%%%%%%%%%%%%%%%%%%%%%%%%

\renewcommand{\thepage}{4-\arabic{page}}
\setcounter{page}{1}
\setcounter{section}{-1}

\section
{
\qquad\qquad\qquad\qquad $\boldsymbol{\S}$\textbf{4}.\quad  Rappel
}
\setlength\parindent{2em}
\setcounter{theoremn}{0}
%%----------------------------------------------------------------------------------------------01
\ \indent

The purpose of this \S \ is to recall the main result of [2-(h)], which will then enable us to put into perspective the central theme of the present paper.  

Thus, let $\alpha$ be a \mK-central, \mK-finite element of $C_c^\infty(G)$ $-$then, subject to the assumptions and conventions of [2-(g), \S 9], 
\[
L_{G / \Gamma}^\dis(\alpha)
\]
is of the trace class, its trace being equal to
\[
\bK(\bH : \alpha : \Gamma)
\]
less
\[
\bFnc(\bH :\bH_\bO : \alpha : \Gamma), 
\]
where (cf. Theorem 4.3 of [2-(h)])
\[
\bFnc(\bH :\bH_\bO : \alpha : \Gamma) 
\]
is equal to
\[
\vsx\sum\limits_{\sC_i, \sC_{i_0}} \ 
\sum\limits_{\bO_{i_0}} \ 
\sum\limits_{w_{i_0}^\dagger}
\]
%2
\[
\vsx\frac
{(-1)^{\dim({}^\prime\mfa_{i_0} (w_{i_0}^\dagger))}}
{(2 \pi)^{\dim(\mfa_{i_0}^\prime (w_{i_0}^\dagger))}}
\]
%3
\[
\vsx\times\ \frac
{1}
{\abs{\det ((1 - w_{i_0}^\dagger) | \Img (1 - w_{i_0}^\dagger ))}}
\]
%4
\[
\vsx\times \ \frac
{1}
{* ({}^\prime \sC_{i_0} (w_{i_0}^\dagger) )\bullet * (\sC_{i_0}^{\dagger\prime} ))}
\]
%5
\[
\vsx\times \ \int_{\sqrt{-1} \ \mfa_{i_0}^\prime (w_{i_0}^\dagger)} \langle 
\bp(\Gamma : \sC_i : \bH - \bH_\bO), 
\]
%6
\[
\vsx\btr_{\sC_i} \bigg( 
\sum\limits_{j_0, {}^\prime w_{j_0 i_0}} \hsx
D_{*}^{w_{j_0}^\dagger} 
[\bc (P_{i_0} | A_{i_0} : P_{j_0} | A_{j_0} : {}^\prime w_{j_0 i_0}^{-1} : ?)
\bullet
\be ( \sC_{j_0} : \bH_\bO : ?) ]_{{}^\prime w_{j_0 i_0} \Lambda_{i_0}^\prime}
\]
%7
\[
\vsx\times \bc ( P_{j_0} | A_{j_0} : P_{i_0} | A_{i_0} : {}^\prime w_{j_0 i_0} : \Lambda_{i_0}^\prime) 
\]
%8
\[
\vsx\times \bc ( P_{i_0} | A_{i_0} : P_{i_0} | A_{i_0} : w_{i_0}^\dagger : 0 )
\bullet 
\bInd_{\sC_{i_0}}^G ( ( \bO_{i_0},\Lambda_{i_0}^\prime)) (\alpha)) \rangle \abs{d \Lambda_{i_0}^\prime}.
\]
Upon setting $\bH_\bO = \bH$ in 
\[
\bFnc(\bH :\bH_\bO : \alpha : \Gamma)
\]
%%----------------------------------------------------------------------------------------------02
we get a polynomial in \bH,
\[
\bCon - \bSp (\bH : \alpha : \Gamma),
\]
that represents the contribution to the trace arising from the continuous spectrum (cf. [2-(h), \S 5]), 

It remains to analyze
\[
\bK(\bH : \alpha : \Gamma).
\]
Referring to [2-(f), p. 1433] for its definition, break up $\Gamma$, 
\[
\Gamma \ = \ \coprod\limits_{\sC_0} \ \coprod\limits_{\sC \succeq \sC_0} \Gamma_{\sC,\sC_0}, 
\]
and write, as is permissible, 
\[
\bK(\bH : \alpha : \Gamma)
\ = \ 
\sum\limits_{\sC_0} \ \sum\limits_{\sC \succeq \sC_0} \bK(\bH : \alpha : \Gamma_{\sC,\sC_0}).
\]
It will therefore be enough to analyze each of the 
\[
 \bK(\bH : \alpha : \Gamma_{\sC,\sC_0})
\]
separately.  
And the easiest of these to handle is the case when $\sC = \sC_0$: 
\[
 \bK(\bH : \alpha : \Gamma_{\sC_0,\sC_0}),
\]
the study of which will occupy us for the remainder of the paper.

%%%%%%%%%%%%%%%%%%%%%%%%%%%%%%%%%%%%%%
%%%%%%%%%%%%%%%%%%%%%%%%%%%%%%%%%%%%%%
%%%%%%%%%%%%%%%%%%%%%%%%%%%%%%%%%%%%%%

\renewcommand{\thepage}{5-\arabic{page}}
\setcounter{page}{1}
\setcounter{section}{-1}

\section
{
\qquad\qquad$\boldsymbol{\S}$\textbf{5}.\quad  Analysis of $\bK(\bH:\balpha:\bGamma_{\sC_0,\sC_0})$
}
\setlength\parindent{2em}
\setcounter{theoremn}{0}
%%----------------------------------------------------------------------------------------------01
\ \indent 
The purpose of this \S\  is to determine the contribution to the trace furnished by 
\[
\bK(\bH : \alpha : \Gamma_{\sC_0,\sC_0}).
\]

As we shall see, the evaluation will ultimately be in terms of weighted orbital integrals 
(or just orbital integrals if $\sC_0 = \{G\}$).

Put $\Delta = \Gamma_{\sC_0,\sC_0}$ (cf. \S3) and, as in [2-(f), \S5], form 
\[
\Phi \hsx = \hsx \{K_{\alpha,\Delta} (P:?) : P \in \sC_\Gamma\}.
\]
Then, by definition, 
\[
\bK(\bH : \alpha : \Delta) 
\hsx = \hsx
\int_{G / \Gamma} Q(\bH : \Phi)(x) d_G(x).
\]
Thanks to what can be found in [2-(f), \S7],
\[
\int_{G / \Gamma} Q(\bH : \Phi)(x) d_G(x)
\]
is a polynomial in \bH.  
Next, as in [2-(f), \S5], form 
\[
\phi \hsx = \hsx \{k_{\alpha,\Delta} (P:?) : P \in \sC_\Gamma\}.
\] 
We shall prove eventually that 
\[
\int_{G / \Gamma} Q(\bH : \phi)(x) d_G(x)
\]
is also a polynomial in \bH if $\bH \ll \textbf{0}$.  
On the other hand, Theorem 5.3 in [2-(f)] implies that the difference
\[
\int_{G / \Gamma} Q(\bH : \Phi)(x) d_G(x) - \int_{G / \Gamma} Q(\bH : \phi)(x) d_G(x)
\]
is o(\bH), hence vanishes.  
So, it will be enough to examine
\[
\int_{G / \Gamma} Q(\bH : \phi)(x) d_G(x).
\]

To this end, keep to the notation of \S3  and decompose $\Delta$ into $\Gamma$-conjugacy classes $\Delta_i$: 
\[
\Delta \hsx = \hsx \coprod \Delta_i.
\]
%%----------------------------------------------------------------------------------------------02
Form anew
\[
\phi_i \hsx = \hsx \{k_{\alpha,\Delta_i} (P:?) : P \in \sC_\Gamma\}.
\]
Then
\[
k_{\alpha,\Delta} \hsx = \hsx \sum\limits_i k_{\alpha,\Delta_i},
\]
a locally finite sum, and, on the basis of Proposition 5.2 in [2-(f)], 
\[
\int_{G / \Gamma} Q(\bH : \phi)(x) d_G(x) 
\hsx = \hsx 
\sum\limits_i \ \int_{G / \Gamma} Q(\bH : \phi_i)(x) d_G(x) ,
\]
thereby reducing our study to that of 
\[
\int_{G / \Gamma} Q(\bH : \phi_0)(x) d_G(x),
\]
$i = 0$ being the fixed index.

An additional reduction is possible provided that we take into account the machinery from \S3.  
Indeed, by definition, at any particular $Q \in \sC_\Gamma$, 
\[
k_{\alpha,\Delta_0}(Q : x) 
\hsx = \hsx 
\sum\limits_{\gamma \in \Delta_0 \cap Q} \ \alpha(x \gamma x^{-1}).
\]
But (cf. \S3),
\[
\Delta_0 \cap Q 
\hsx = \hsx 
\coprod\limits_{P \in E(\gamma_0)} \ \Delta_0(P)Q.
\]
The cardinality of the $\Gamma_{\gamma_0}$-equivalence class of \mP in $\sP(L_0)$ is 
\[
[\Gamma_{\gamma_0} : \Gamma_{\gamma_0} \cap P].
\]
It therefore follows that
\[
k_{\alpha,\Delta_0}(Q : x) 
\hsx = \hsx 
\sum\limits_{P \in \sP(L_0)} \ 
\bigg( \frac{1}{[\Gamma_{\gamma_0} : \Gamma_{\gamma_0} \cap P]} 
\ \bullet 
\sum\limits_{\gamma \in \Delta_0(P)_Q} \ \alpha(x \gamma x^{-1}) \bigg).
\]
In this connection, bear in mind that 
\[
\Delta_0(P)_Q \hsx \neq \hsx \emptyset
\]
iff \mP and \mQ are $\Gamma$-conjugate, so we have not really overloaded the sum.  
Each summand determines a $\Gamma$-compatible family of functions on \mG (cf. Lemma 3.1), call it $\Phi(P)$.  
The focal point of the analysis thus becomes
\[
\int_{G /\Gamma} \ \sum\limits_{P \in \sP(L_0)} \ Q(\bH : \Phi(P))(x) d_G(x).
\]
%%----------------------------------------------------------------------------------------------03

Individually, 
\[
\int_{G /\Gamma} Q(\bH : \Phi(P))(x) d_G(x)
\]
may well diverge.  
To deal with this difficulty, we shall need a lemma.

\vspace{0.3cm}

{\small\bf Lemma 5.1} \ %5.1  
Let $f$ be a bounded, compactly supported, measurable function on $G /  \Gamma$ $-$then
\[
\int_{G /\Gamma} Q(\bH : \Phi(P))(x) f(x) d_G(x)
\]
is equal to
\[
=\ 
\frac{(-1)^{\rank(P)}}{[\Gamma_{\gamma_0} : \Gamma_{\gamma_0} \cap P]}
\ \bullet
\int_{G /\Gamma_{\gamma_0} \cap L_0} \ \alpha(x \gamma_0 x^{-1}) f(x) 
\]
\[
\times \chi_{P,A}: \rsC (\bH(P) - H_{\restr{P}{A}}(x)) d_G(x).
\]

[Note: \ Here, the split component \mA of \mP is per

\[
(P,S;A) \hsx \succeq (P_0,S_0;A_0) \qquad (P_0 \in \sP_0(L_0)),
\]
where $L_0 = M_0 \bullet A_0$ with $\gamma_0 \in M)o$ (cf. \S3).]

{\small\bf Proof.}  
In fact, 

\[
\int_{G /\Gamma} \  Q(\bH : \Phi(P)) (x) f(x) d_G(x)
\]
%2
\[
=\ 
\frac{(-1)^{\rank(P)}}{[\Gamma_{\gamma_0} : \Gamma_{\gamma_0} \cap P]}
\ \bullet
\int_{G /\Gamma} \ 
\sum_{\delta \in \Gamma / \Gamma \cap P} \ 
\sum_{\gamma \in \Delta_0(P)_P} 
\]
\[
\times \ \alpha(x \delta \gamma  \delta^{-1} x^{-1}) f(x) \chi_{P,A}: \rsC (\bH(P) - H_{\restr{P}{A}}(x \gamma)) d_G(x)%done
\]
%3
\[
=\ 
\frac{(-1)^{\rank(P)}}{[\Gamma_{\gamma_0} : \Gamma_{\gamma_0} \cap P]}
\ \bullet
\int_{G /\Gamma \cap P} \ \sum_{\gamma \in \Delta_0(P)_P} 
\]
\[
\times \ \alpha(x \gamma  x^{-1}) f(x) \chi_{P,A}: \rsC (\bH(P) - H_{\restr{P}{A}}(x)) d_G(x) %done
\]
%4
\[
=\ 
\frac{(-1)^{\rank(P)}}{[\Gamma_{\gamma_0} : \Gamma_{\gamma_0} \cap P]}
\ \bullet
\int_{G /\Gamma \cap P} \ 
\sum\limits_{\gamma \in \Gamma \cap P /  \Gamma_{\gamma_0} \cap P}
\]
\[
\times \ \alpha(x \gamma \gamma_0 \gamma^{-1} x^{-1}) f(x) \chi_{P,A}: \rsC (\bH(P) - H_{\restr{P}{A}}(x)) d_G(x) %d
\]
%5
\[
=\ 
\frac{(-1)^{\rank(P)}}{[\Gamma_{\gamma_0} : \Gamma_{\gamma_0} \cap P]}
\ \bullet
\int_{G /\Gamma_{\gamma_0} \cap P} \ \alpha(x \gamma_0 x^{-1}) f(x) %done
\]
\[
\times \chi_{P,A}: \rsC (\bH(P) - H_{\restr{P}{A}}(x)) d_G(x) %done
\]
%%----------------------------------------------------------------------------------------------04
%6
\[
=\ 
\frac{(-1)^{\rank(P)}}{[\Gamma_{\gamma_0} : \Gamma_{\gamma_0} \cap L_0]}
\ \bullet
\int_{G /\Gamma_{\gamma_0} \cap L_0} \ \alpha(x \gamma_0 x^{-1}) f(x) %done
\]
\[
\times \chi_{P,A}: \rsC (\bH(P) - H_{\restr{P}{A}}(x)) d_G(x), %done
\]
as desired.

\vspace{0.3cm}

If $\bH \ll \textbf{0}$, then 
\[
Q(\bH : \phi_0)
\]
has bounded support (cf. Proposition 5.2 in [2-(f)]), so, if $f \equiv 1$ on a large enough set, then 
\[
\int_{G /\Gamma} Q(\bH : \phi_0)(x) d_G(x)
\]
is equal to
\[
\frac{1}{[\Gamma_{\gamma_0} : \Gamma_{\gamma_0} \cap L_0]}
\bullet
\int_{G /\Gamma_{\gamma_0} \cap L_0} \ \alpha(x \gamma_0 x^{-1}) f(x)
\]
%2
\[
\times 
\bigg ( \sum\limits_{P \in \sP(L_0)} \ (-1)^{\rank(P)} \chi_{P,A} \rsC (\bH(P) - H_{\restr{P}{A}}(x)) \bigg ) d_G(x)
\]
or still 
\[
\frac{1}{[\Gamma_{\gamma_0} : \Gamma_{\gamma_0} \cap L_0]}
\bullet
\int_{G / A_0 \bullet (\Gamma_{\gamma_0} \cap L_0) } \ \alpha(x \gamma_0 x^{-1}) \ \int\limits_{\mfa_0 }f(xe^H)
\]
%2
\[
\times 
\bigg ( \sum\limits_{P \in \sP(L_0)} \ (-1)^{\rank(P)} \chi_{P,A} \rsC (\bH(P) - H_{\restr{P}{A}}(x) - H) \bigg ) dH d_G(x).
\]
We shall see in a bit how the ``$f$'' can be eliminated.  
Anticipating this, let us consider in more detail 
\[
\int\limits_{\mfa_0 } 
\big ( \sum\limits_{P \in \sP(L_0)}  \ (-1)^{\rank(P)} 
\chi_{P,A} \reflectbox{$\sC$} (\bH(P) - H_{\restr{P}{A}}(x) - H) \big ) dH.
\]

As might be expected, the issue is primarily combinatorial in character.  
That being the case, our basic tools will be drawn from the repository in [2-(b), \S2] and [2-(d), \S3], 
the notation of which will be employed below without further comment.

It is well-known and familiar that there is a one-to-one correspondence between
\begin{align*}
&(i) \quad\ \  \{P \in \sP(L_0)\};\\
&(ii) \quad \ \{\mfa, \sC_\mfa\};\\
&(iii) \quad \{W, \sC_W\}.
\end{align*}

%%----------------------------------------------------------------------------------------------05
The following lemma gives yet another parametrization.  
Put $\Phi = \sum (\mfg, \mfa_0)$.

\vspace{0.3cm}

{\small\bf Lemma 5.2} \ % 5.2
Fix a chamber $\sC_\bO$ in $\mfa_0$ $-$then the set of all pairs $(W,\sC_W)$ is in a one-to-one correspondence with the 
set of all pairs $(\sC_0,F_0)$, where 
\[
F_0 \subset F_0(\sC_0) \cap \Phi^+(\sC_0)
\]
and 
\[
\begin{cases}
\ W = \mfa_0(F_0)\\
\ (\sC(F_0),\sC_W) \longleftrightarrow \sC_0
\end{cases}
.
\]

\vspace{0.1cm}

{\small\bf Proof.}  
It suffices to show that for each $W$ there exists a unique $\sC(W)$ such that $F_0(\sC(W)) \subset \Phi^+(\sC_0)$.  
But the projection of $\sC_\bO$ onto $W$ is a connected open set on which an element of $\Phi(W)$ is either positive 
or negative, hence is contained in a unique chamber $\sC(W)$.
And $\sC(W)$ will do.  
If ?$(W)$ is another chamber, some member of $F_0(?(W))$ must be negative on $\sC(W)$, thus negative on the projection of 
$\sC_\bO$ and so negative on $\sC_\bO$.  
This means that $\sC(W)$ is the only chamber that will do.

\vspace{0.1cm}

In the terminology of Arthur, if $\bH \ll \textbf{0}$, then 
\[
\{\bH(P_0) - H_{\restr{P_0}{A_0}}(x) : P_0 \in \sP_0(L_0)\}
\]
is a negative $A_0$-orthogonal set (cf. [1-(a), p. 221]).  
The elements $\sC_0$ of $\sC_0(A_0)$  are in a one-to-one correspondence with the elements $P_0$ of $\sP_0(L_0)$ 
(cf. TES, p. 66).  
Assuming that $\sC_0 \longleftrightarrow P_0$, put 
\[
T_{\sC_0}(x : \bH) \hsx = \hsx \bH(P_0) - H_{\restr{P_0}{A_0}}(x). 
\]
 
We can now get on with the manipulation.  
Thus
%1
\[
\sum\limits_{P \in \sP(L_0)}
(-1)^{\rank(P)}
\chi_{P,A :} \rsC (\bH(P) - H_{\restr{P_0}{A_0}}(x) - H)
\]
%2
\[
=\ 
\sum\limits_{\mfa}
\sum\limits_{\sC_\mfa}
(-1)^{\rank(P)}  \chi^{*,F_0(\sC_\mfa)}
(T_{\sC_0} (x : \bH) - H)
\]
%3
\[
=\ 
\sum\limits_{W, \sC_W} \ 
(-1)^{\ell_0 - \dim(W)} \ \chi^{F_0(\sC(W)), F_0(\sC_0)} (T_{\sC_0} (x : \bH) - H)
\]
%4
\[
=\ 
\sum\limits_{\sC_0}
\sum\limits_{\{F:F \subset F_0(\sC_0) \cap \Phi^+(\sC_0)\}} \ 
(-1)^{\ell_0 -\#(F)}   \chi^{F, F_0(\sC_0)} (T_{\sC_0} (x : \bH) - H).
\]
The sum over \mF can be cut down considerably.  
Indeed, if \mX is a variable, then 
\[
\chi^{F, F_0(\sC_0)}(X) \hsx = \hsx
\begin{cases}
\ 1 \qquad \text{if} \ \lambda^i(X) > 0 \ \forall \ \lambda_i \in F_0(\sC_0) - F\\
\ 0 \qquad \text{otherwise}
\end{cases}
.
\]
%%----------------------------------------------------------------------------------------------06
So, with 
\[
F_0(X) 
\hsx = \hsx 
\{\lambda_i \in F_0(\sC_0) : \lambda^i(X) \leq 0\},
\]
%2
\[
\sum\limits_{\{F: F \subset F_0(\sC_0) \cap \Phi^+(\sC_0)\}} \ 
(-1)^{\ell_0 - \#(F)} \chi^{F,F_0(\sC_0)}(X)
\]
%3
\[
= \ 
\sum\limits_{\{F: F_0(X) \subset F \subset F_0(\sC_0) \cap \Phi^+(\sC_0)\}} \ 
(-1)^{\ell_0 - \#(F)} 
\]
%4
\[
=\  
\begin{cases}
\ (-1)^{\ell_0 - \#(F_0(\sC_0) \cap \Phi^+(\sC_O))} \qquad \text{if} \ F_0(X) = F_0(\sC_0) \cap \Phi^+(\sC_O)\\
\ 0 \qquad\qquad\qquad\qquad\qquad\qquad \text{otherwise}
\end{cases}
\]
or still
\[
\sum\limits_{\{F: F \subset F_0(\sC_0) \cap \Phi^+(\sC_0)\}} \ 
(-1)^{\ell_0 - \#(F)} \chi^{F,F_0(\sC_0)}(X)
\]
%2
\[
=\ (-1)^{\ell_0 - \#(F_0(\sC_0) \cap \Phi^+(\sC_O))}
\]
%3
\[
\times \ \tau_{*,F_0(\sC_0)} (F_0(\sC_0) - \Phi^+(\sC_O) : X).  %%come back and fix C_0 and C_O
\]
Inserting this then leads to 
\[
\sum\limits_{\sC_0} \ (-1)^{\#(F_0(\sC_0) - \Phi^+(\sC_O))}
\]
%2
\[
\times \ \tau_{*,F_0(\sC_0)} (F_0 (\sC_0) - \Phi^+ (\sC_O) : T_{\sC_0} (x : \bH) - H),
\]
the characteristic function of the convex hull of the $T_{\sC_0}(x : \bH)$ (cf. [1-(a), pp. 218-219]).  
Integrating over $\mfa_0$ then gives its volume, call it
\[
v_{\mfa_0} (x : \bH).
\]

\vspace{0.3cm}

In summary, therefore, if we ignore the ``$f$'', then 
\[
\int\limits_{G / \Gamma} \ Q(\bH : \phi_0)(x) d_G(x)
\]
is equal to 
\[
\frac{1}{\Gamma_{\gamma_0} : \Gamma_{\gamma_0} \cap L_0]} 
\ \bullet 
\int_{G / A_0 \bullet (\Gamma_{\gamma_0} \cap L_0)} \ 
\alpha(x \gamma_0 x^{-1}) v_{\mfa_0} (x : \bH) d_G(x).
\]
And:

\vspace{0.3cm}

{\small\bf Lemma 5.3} \ %5.3
The integral
\[
\int_{G / A_0 \bullet(\Gamma_{\gamma_0} \cap L_0)} \ 
\alpha (x \gamma_0 x^{-1}) v_{\mfa_0} (x : \bH) d_G(x)
\]
%%----------------------------------------------------------------------------------------------07
is compactly supported.

\vspace{0.1cm}

{\small\bf Proof.}
The centralizer of $\gamma_0$ in $M_0$ is the same as the centralizer of $\gamma_0$ in $S_0$.  
Moreover, $\gamma_0$ is $M_0$-regular, thus the quotient
\[
(M_0) _{\gamma_0} / \Gamma_{\gamma_0} \cap M_0
\]
is compact (cf. \S 2: $\Gamma \cap M_0$) has finite index in $\Gamma_{M_0}$).  
So, 
\[
\int_{G / A_0 \bullet(\Gamma_{\gamma_0} \cap L_0)} \ 
\alpha (x \gamma_0 x^{-1}) v_{\mfa_0} (x : \bH) d_G(x) 
\]
\begin{align*}
\vsx&= \ \int_{S_0 / \Gamma_{\gamma_0} \cap M_0} \ 
\alpha (s \gamma_0 s^{-1}) v_{\mfa_0} (s : \bH) d_S(s) \\
&= \vol((M_0)_{\gamma_0} / \Gamma_{\gamma_0} \cap M_0)
\end{align*}
\[
\times \ \int_{S_0 / (S_0)_{\gamma_0}}
\alpha (s \gamma_0 s^{-1}) v_{\mfa_0} (s : \bH) d_S(s).
\]
It remains only to note that the $S_0$-orbit of $\gamma_0$ is the $M_0$-orbit of $\gamma_0$ times $N_0$, which is closed.

\vspace{0.3cm}

Consequently, the ``$f$'' can in fact be dispensed with.

\vspace{0.3cm}

{\small\bf Remark.}  
One would like to say that 
\[
\int_{G / A_0 \bullet(\Gamma_{\gamma_0} \cap L_0)} \ 
\alpha (x \gamma_0 x^{-1}) v_{\mfa_0} (x : \bH) d_G(x)
\]
is a weighted orbital integral.  This is certainly the case if $G_{\gamma_0}$ is contained in $L_0$.  
But, as we have already noted in \S2, this need ot be true in general although $G_{\gamma_0}^0$ 
 is always contained in $L_0$ (cf. Lemma 2.6).  
On the other hand, it is not difficult to see that
\[
G_{\gamma_0} / A_0 \bullet (\Gamma_{\gamma_0} \cap L_0)
\]
is at least compact, thus
\[
\int_{G / A_0 \bullet(\Gamma_{\gamma_0} \cap L_0)} \ 
\alpha (x \gamma_0 x^{-1}) v_{\mfa_0} (x : \bH) d_G(x)
\]
\[
= \ 
\int_{G / G_{\gamma_0}} \ 
\alpha (x \gamma_0 x^{-1}) 
A \ver_{\gamma_0}(v_{\mfa_0} (x : \bH)) d_G(x),
\]
%%----------------------------------------------------------------------------------------------08
the density 
\[
A \ver_{\gamma_0}(v_{\mfa_0} (x : \bH)) 
\hsx = \hsx 
\int_{G_{\gamma_0} / A_0 \bullet(\Gamma_{\gamma_0} \cap L_0)} \ 
v_{\mfa_0} (xy : \bH) d_{G_{\gamma_0}}(y)
\]
being an averaged volume element.

\vspace{0.3cm}

Having evaluated 
\[
\int_{G / \Gamma} \ Q(\bH : \phi_0)(x) d_G(x)
\]
in closed form, the next step in the analysis is to prove that it is a polynomial in \bH if $\bH \ll \textbf{0}$

For this purpose, recall that from the 
\[
\begin{cases}
\ -H_{\restr{P_0}{A_0}}(x)\\
\ \bH(P_0) - H_{\restr{P_0}{A_0}} (x)
\end{cases}
\qquad (P_0 \in \sP_0(L_0))
\]
one can manufacture Detroit families
\[
\begin{cases}
\ \be^{-H_?(x)}\\
\ \be^\bH(?) - H_?(x)
\end{cases}
\]
by exponentiation (cf. [2-(d), p. 163]).  
Put
\[
\begin{cases}
\ v_{\mfa_0} (x) \ = \ (-1)^{\ell_0}\Sh_{\be^{- H_?(x)}} (0)\\
\ v_{\mfa_0} (x : \bH) \ = \ (-1)^{\ell_0}\Sh_{\be^\bH(?) - H_?(x)} (0)
\end{cases}
,
\]
a permissible agreement.  
Owing now to Corollary 3.2 in [2-(h)], we have
\begin{align*}
v_{\mfa_0}(x : \bH) \hsx \ 
\vsy&=\  \ (-1)^{\ell_0}\Sh_{\be^\bH(?) - H_?(x)} (0)\\
\vsy&= \ (-1)^{\ell_0} \ \sum\limits_{W,\sC_W} \ \Sh_{\be^{-H_?(x)}(\sC_W)} (0) \bullet A_{\be^\bH,\sC_W}(0)\\
\vsy&= \ \sum\limits_{P \in \sP(L_0)} v_{\mfa_0^\dagger} (m_x) \bullet p(\Gamma : P : \bH(P)),
\end{align*}
where $v_{\mfa_0^\dagger}$ is the daggered analogue of $v_{\mfa_0}$ and $p(\Gamma : P : ?)$ 
is an Arthur polynomial (cf. [2-(f), p. 1429]).  Accordingly, 
\[
\int_{G / A_0 \bullet (\Gamma_{\gamma_0} \cap L_0)} \ 
\alpha (x \gamma_0 x^{-1}) v_{\mfa_0} (x : \bH) d_G(x)
\]
%%----------------------------------------------------------------------------------------------09
is the sum over the $P \in \sP(L_0)$ of the
\[
\int_{G / A_0 \bullet (\Gamma_{\gamma_0} \cap L_0)} \ 
\alpha (x \gamma_0 x^{-1}) v_{\mfa_0^\dagger} (m_x)d_G(x)
\]
times 
\[
p(\Gamma : P : \bH(P)).
\]
And so
\[
\int_{G / \Gamma} Q(\bH : \phi_0) (x) d_G(x)
\]
is in fact a polynomial in \bH if $\bH \ll \textbf{0}$.

While we are at it, let us also observe that
\[
\frac{1}{[\Gamma_{\gamma_0} : \Gamma_{\gamma_0} \cap L_0]}
\ \bullet
\int_{G / A_0 \bullet (\Gamma_{\gamma_0} \cap L_0)} \ 
\alpha (x \gamma_0 x^{-1}) v_{\mfa_0^\dagger} (m_x)d_G(x)
\]
\[
= \ 
\frac{1}{[\Gamma_{\gamma_0} : \Gamma_{\gamma_0} \cap L_0]}
\ \bullet
\int_N 
\int_{M / A_0^\dagger (\Gamma_{\gamma_0} \cap L_0)} 
\alpha(n m \gamma_0 m^{-1}n^{-1}) v_{\mfa_0^\dagger} d_M(m) d_N(n)
\]
\[
= \ 
\frac{1}{[\Gamma_{\gamma_0} : \Gamma_{\gamma_0} \cap M]}
\bullet
\frac{1}{\abs{\det(\Ad(\gamma_0) | \bn - 1)}}
\]
\[
\times \ 
\int_{M / A_0^\dagger \bullet (\Gamma_{\gamma_0} \cap M)} 
\alpha^P (m \gamma_0 m^{-1}) v_{\mfa_0^\dagger} (m) d_M(m), 
\] 
the last integral being the \mM-analogue of the integral
\[
\int_{G / A_0 \bullet \Gamma_{\gamma_0}} \ 
\alpha (x \gamma_0 x^{-1}) v_{\mfa_0} (x)d_G(x)
\]
on \mG.

\vspace{0.3cm}

{\small\bf Remark.}
When $\mfa_0$ is special, one can interpret $v_{\mfa_0}(x)$ geometrically in that it then gives the volume of the convex hull 
of the $\{-H_{\restr{P_0}{A_0}}(x)\}$ but this will not be true in general.

Our initial objective can now be realized.  
Because the support of $\alpha$ is compact, Proposition 3.3 guarantees us that only finitely many
\[
\{\gamma\}_\Gamma : \alpha \in \Gamma_{\sC_0,\sC_0}
\]
actually intervene, so all $\bH \ll \textbf{0}$ will work for each of them simultaneously.  
Hence:

\vspace{0.1cm}

%%----------------------------------------------------------------------------------------------10
{\small\bf Theorem 5.4} \ %5.4
$\bK(\bH : \alpha : \Gamma_{\sC_0,\sC_0})$ is equal to
\[
\sum\limits_{\{\gamma\}_\Gamma : \gamma \in \Gamma_{\sC_0,\sC_0}}
\frac{1}{[\Gamma_{\gamma} : \Gamma_{\gamma} \cap L_0(\gamma)]}
\]
\[
\times \ 
\int_{G / A_0(\gamma) \bullet (\Gamma_\gamma \cap L_0(\gamma))} 
\alpha (x \gamma x^{-1}) v_{\mfa_0(\gamma)} (x : \bH) d_G(x)
\]
or still
\[
\sum\limits_{\{\gamma\}_\Gamma : \gamma \in \Gamma_{\sC_0,\sC_0}}
\sum\limits_{P \in \sP(L_0(\gamma))} 
\frac{1}{[\Gamma_{\gamma} : \Gamma_{\gamma} \cap M]}
\bullet
\frac{1}{\abs{\det(\Ad(\gamma) | \bn - 1)}}
\]
\[
\times p (\Gamma : P : \bH(P))
\]
\[
\times \ 
\int_{M / A_0^\dagger(\gamma) \bullet (\Gamma_{\gamma} \cap M)} 
\alpha^P (m \gamma m^{-1}) v_{\mfa_0^\dagger(\gamma)} (m) d_M(m).
\]

[Note: \ For the sake of brevity, we have written $P = M \bullet A \bullet N$ in place of 
$P(\gamma) = M(\gamma) \bullet A(\gamma) \bullet N(\gamma)$.]

%%%%%%%%%%%%%%%%%%%%%%%%%%%%%%%%%%%%%%
%%%%%%%%%%%%%%%%%%%%%%%%%%%%%%%%%%%%%%
%%%%%%%%%%%%%%%%%%%%%%%%%%%%%%%%%%%%%%

\renewcommand{\thepage}{6-\arabic{page}}
\setcounter{page}{1}
\setcounter{section}{-1}

\section
{
\qquad\qquad$\boldsymbol{\S}$\textbf{6}.\quad  Passage to Standard Form
}
\setlength\parindent{2em}
\setcounter{theoremn}{0}
%%----------------------------------------------------------------------------------------------01
\ \indent 
The purpose of this \S\  is to recast the expression obtained in Theorem 5.4 for
\[
\bK(\bH:\alpha:\Gamma_{\sC_0,\sC_0})
\]
so as to reflect the presence of the 
$P_{i \mu} = M_{i \mu} \bullet A_{i \mu} \bullet N_{i \mu}$ 
with 
$A_{i \mu}$ special (cf [2-(a), pp. 65-70]), the point being that these are the parabolics of reference.  
This issue is therefore primarily one of bookkeeping, albeit a little on the involved side.

We can evidently write
\[
\bK(\bH:\alpha:\Gamma_{\sC_0,\sC_0})
\]
in the form 
\[
\vsx\sum\limits_{\{\gamma\}_\Gamma : \gamma \in \Gamma_{\sC_0,\sC_0}}
\ 
\sum\limits_{P \in \sP(L_0(\gamma))} 
\frac{1}{[\Gamma_\gamma : \Gamma_\gamma \cap M_0(\gamma)]}
\]
\[
\times p(\Gamma : P : \bH(P)) 
\]
\[
\times \int_{G / A_0(\gamma) \bullet (\Gamma_\gamma \cap M_0(\gamma))} 
\alpha(x \gamma x^{-1}) v_{\mfa_0^\dagger(\gamma)} (x) d_G(x),
\]
the effect of which is to base part of the data at the $L_0(\gamma)$-level.  
Now introduce a parameter $\bH_\bO \in \mfa$ $-$then still another way to write \\
\[
\bK(\bH:\alpha:\Gamma_{\sC_0,\sC_0})
\]
is 
\[
\vsx\sum\limits_{\{\gamma\}_\Gamma : \gamma \in \Gamma_{\sC_0,\sC_0}}
\ 
\sum\limits_{P \in \sP(L_0(\gamma))} 
\frac{1}{[\Gamma_\gamma : \Gamma_\gamma \cap M_0(\gamma)]} 
\]
\[
\times p(\Gamma : P : \bH(P) - \bH_\bO(P)) 
\]
\[
\times \int_{G / A_0(\gamma) \bullet (\Gamma_\gamma \cap M_0(\gamma))} 
\alpha(x \gamma x^{-1}) v_{\mfa_0^\dagger(\gamma)} (x : I_M(\bH_\bO)) d_G(x).
\]
As we shall see, the rationale here is that the introduction of $\bH_\bO$ allows certain transitions to take place without the need for the introduction of compensating factors.

Consider an arbitrary $P \in \sP(L_0(\gamma))$.  
Choose $\gamma(i:\mu) \in \Gamma$ such that 
\[
\gamma(i:\mu) P \gamma(i:\mu)^{-1} \hsx = \hsx P_{i \mu}
\]
and let
\[
\gamma_{i \mu} \hsx = \hsx \gamma(i:\mu)  \gamma \gamma(i:\mu)^{-1}.
\]
%%----------------------------------------------------------------------------------------------02
Obviously, 
\[
P_{i \mu} \in \sP(L_0(\gamma_{i \mu}))
\]
but there is no reason to expect that 
\[
 L(\gamma_{i \mu}) \hsx = \hsx M(\gamma_{i \mu}) \bullet A(\gamma_{i \mu})
\]
is the special $L_{i \mu} = M_{i \mu} \bullet A_{i \mu}$.

It is clear that 
\[
p(\Gamma : P : \bH(P) - \bH_\bO(P)) 
\hsx = \hsx
p(\Gamma : P_{i \mu} : \bH(P_{i \mu}) - \bH_\bO(P_{i \mu})).
\]
Additionally, one can verify that
\[
\int_{G / A_0(\gamma) \bullet (\Gamma_\gamma \cap M_0(\gamma))} 
\alpha(x \gamma x^{-1}) v_{\mfa_0^\dagger(\gamma)} (x : I_M(\bH_\bO)) d_G(x)
\]
\[
= \ \int_{G / A_0(\gamma_{i \mu}) \bullet (\Gamma_{\gamma_{i \mu}} \cap M_0(\gamma_{i \mu}))} 
\alpha(x \gamma_{i \mu} x^{-1}) v_{\mfa_0^\dagger(\gamma_{i \mu})} 
(x : I_{M(_{i \mu})}(\bH_\bO)) d_G(x).
\]
These are the two crucial relations which allow the passage from $\gamma$ to $\gamma_{i \mu}$.  
They would fail to hold if not for the presence of $\bH_\bO$.

Accordingly, integrating out $N_{i \mu}$, 
\[
\bK(\bH:\alpha:\Gamma_{\sC_0,\sC_0})
\]
becomes
\[
\sum\limits_{\{\gamma\}_\Gamma : \gamma \in \Gamma_{\sC_0,\sC_0}}
\ 
\sum\limits_{P \in \sP(L_0(\gamma))} 
\frac{1}{[\Gamma_{\gamma_{i \mu}} : \Gamma_{\gamma_{i \mu}} \cap M_0(\gamma_{i \mu})]} %\\
\]
%\[
%\vsx\times \frac{1}{\abs{\det(\Ad(\gamma_{i \mu})| \mfn_{i \mu} - 1)}} 
%\]
\[
\vsx\times \frac{1}{\abs{\det(\restr{\Ad(\gamma_{i \mu})}{\mfn_{i \mu}} - 1)}} 
\]
\[
\times p(\Gamma : P_{i \mu} : \bH(P_{i \mu}) - \bH_\bO(P_{i \mu})) %\\
\]
\[
\times \int_{M(\gamma_{i \mu}) / A_0^\dagger(\gamma_{i \mu}) 
\bullet (\Gamma_{\gamma_{i \mu}} \cap M_0(\gamma_{i \mu}))}
\alpha^{P_{i \mu}}(m\gamma_{i \mu} m^{-1}) v_{\mfa_0^\dagger(\gamma_{i \mu})} 
(m : I_{M(\gamma_{i \mu})} (\bH_\bO)) d_{M(\gamma_{i \mu})} (m).
\]

Next, determine $n_{i \mu} \in N_{i \mu}$ so that 
\[
n_{i \mu} A (\gamma_{i \mu}) n_{i \mu}^{-1} \hsx = \hsx A_{i \mu}.
\]
Passing now from $\gamma_{i \mu}$ to 
\[
\delta_{i \mu} \hsx = \hsx n_{i \mu} \gamma_{i \mu} n_{i \mu}^{-1} \in \Gamma_{M_{i \mu}},
\]
%%----------------------------------------------------------------------------------------------03
and agreeing to write
\[
\begin{cases}
\ M_0(\delta_{i \mu}) \quad \text{instead of} \  n_{i \mu} M_0(\gamma_{i \mu}) n_{i \mu}^{-1}\\
\ A_0(\delta_{i \mu}) \quad \text{instead of} \  n_{i \mu} A_0(\gamma_{i \mu}) n_{i \mu}^{-1}
\end{cases}
,
\]
the foregoing then reads
\[
\sum\limits_{\{\gamma\}_\Gamma : \gamma \in \Gamma_{\sC_0,\sC_0}}
\sum\limits_{P \in \sP(L_0(\gamma))} 
\frac
{
[(\Gamma_{M_{i \mu}})_{\delta_{i \mu}} \cap M_0(\delta_{i \mu}) : 
n_{i \mu} \Gamma_{\gamma_{i \mu}} n_{i \mu}^{-1} \cap M_0(\delta_{i \mu})]
}
{
[n_{i \mu} \Gamma_{\gamma_{i \mu}} n_{i \mu}^{-1}  : 
n_{i \mu} \Gamma_{\gamma_{i \mu}} n_{i \mu}^{-1} \cap M_0(\delta_{i \mu})]
}
\]
\[
\vsx\times \frac{1}{\abs{\det(\restr{\Ad(\gamma_{i \mu})}{\mfn_{i \mu}} - 1)}} 
\]
\[
\times p(\Gamma : P_{i \mu} : \bH(P_{i \mu}) - \bH_\bO(P_{i \mu})) %\\
\]
\[
\times \int_{M(\gamma_{i \mu}) / A_0^\dagger(\delta_{i \mu}) 
\bullet ((\Gamma_{M_{i \mu}})_{\delta_{i \mu}} \cap M_0(\delta_{i \mu}))}
\alpha^{P_{i \mu}}(m\delta_{i \mu} m^{-1}) v_{\mfa_0^\dagger(\delta_{i \mu})} 
(m : I_{M_{i \mu}} (\bH_\bO)) d_{M_{i \mu}} (m).
\]

To get an inductive object out of this expression, it will be necessary to replace
\[
\sum\limits_{\{\gamma\}_\Gamma : \gamma \in \Gamma_{\sC_0,\sC_0}}
\ 
\sum\limits_{P \in \sP(L_0(\gamma))} 
\]
by 
\[
\sum\limits_{i, \mu} \ \sum\limits_{\{\delta_{i \mu}\}\Gamma_{M_{i \mu}}}.
\]
As a preliminary, we remark that our notation is slightly deceptive in that $\delta_{i \mu}$ really depends on $\gamma$.  
However, a choice of 

\[
\begin{cases}
\ \gamma \in \Gamma_{\sC_0,\sC_0}\\
\ P \in \sP(L_0(\gamma))
\end{cases}
\]
singles out uniquely
\[
P_{i \mu} \quad \text{and} \quad \{\delta_{i \mu}\}_{\Gamma_{M_{i \mu}}}.
\]
This said, recall that via the daggering procedure $\sC_0$ will induce a disjoint union
\[
\sC_0^\dagger \  \equiv \  \coprod \ \sC_0^\dagger (i : \mu).
\]
Of course, 
\[
\{\delta_{i \mu}\}_{\Gamma_{M_{i \mu}}} \in (\Gamma_{M_{i \mu}})_{\sC_0^\dagger,\sC_0^\dagger} 
\  \equiv \ 
\coprod \ (\Gamma_{M_{i \mu}})_{\sC_0^\dagger(i : \mu),\sC_0^\dagger(i : \mu)}.
\]
And, the ambient 
$\Gamma_{M_{i \mu}}$-Levi subgroup of the ambient 
$\Gamma_{M_{i \mu}}$-cuspidal parabolic in $\sP_0^\dagger(L_0(\delta_{i \mu}))$ is 
$A_0^\dagger(\delta_{i \mu}) \bullet M_0(\delta_{i \mu})$.

\vspace{0.2cm}
%%----------------------------------------------------------------------------------------------04

Consequently, 
\[
\bK(\bH:\alpha:\Gamma_{\sC_0,\sC_0})
\]
reduces to
\[
\sum\limits_{i, \mu}
\ 
\sum\limits_{\{\delta_{i \mu}\}_{\Gamma_{M_{i \mu}}} : 
\delta_{i \mu} \in (\Gamma_{M_{i \mu}})_{\sC_0^\dagger,\sC_0^\dagger}}
C(P_{i \mu}, \{\delta_{i \mu}\}_{\Gamma_{M_{i \mu}}})
\]
\[
\times p(\Gamma : P_{i \mu} : \bH(P_{i \mu}) - \bH_\bO(P_{i \mu})) 
\]
\[
\times \int_{M_{i \mu} / A_0^\dagger(\delta_{i \mu}) 
\bullet 
((\Gamma_{M_{i \mu}})_{\delta_{i \mu}} \cap M_0(\delta_{i \mu}))}
\alpha^{P_{i \mu}}(m\delta_{i \mu} m^{-1}) v_{\mfa_0^\dagger(\delta_{i \mu})} 
(m : I_{M_{i \mu}} (\bH_\bO)) d_{M_{i \mu}} (m).
\]
$\vsx$Here, the constant
\[
C(P_{i \mu}, \{\delta_{i \mu}\}_{\Gamma_{M_{i \mu}}})
\]
is, by definition, the sum
\[
\sum \frac
{
[(\Gamma_{M_{i \mu}})_{\delta_{i \mu}} \cap M_0(\delta_{i \mu}) : 
n_{i \mu} \Gamma_{\gamma_{i \mu}} n_{i \mu}^{-1} \cap M_0(\delta_{i \mu})]
}
{
[n_{i \mu} \Gamma_{\gamma_{i \mu}} n_{i \mu}^{-1}  : 
n_{i \mu} \Gamma_{\gamma_{i \mu}} n_{i \mu}^{-1} \cap M_0(\delta_{i \mu})]
}
\]
\[
\times \ \frac{1}{\abs{\det(\Ad(\gamma_{i \mu}) | \mfn_{i \mu} -1}}
\]
over all possible 
\[
\begin{cases}
\ \gamma \in \Gamma_{\sC_0,\sC_0}\\
\ P \in \sP(L_0(\gamma))
\end{cases}
\]
which gives rise to a fixed
\[
P_{i \mu} \quad \text{and} \quad \{\delta_{i \mu}\}_{\Gamma_{M_{i \mu}}}.
\]
The main technical claim is then: 

\vspace{0.3cm}

{\small\bf Lemma 6.1} \ %01
$C(P_{i \mu},\{\delta_{i \mu}\}_{\Gamma_{M_{i \mu}}})$ is equal to 
\[
\frac{1}
{[(\Gamma_{M_{i \mu}})_{\delta_{i \mu}} : (\Gamma_{M_{i \mu}})_{\delta_{i \mu}} \cap M_0(\delta_{i \mu})]}
.
\]

\vspace{0.1cm}

{\small\bf Proof.} \ 
The way to keep track of what's coming and going is to look at the nonempty 
\[
\Delta_0(P)_{P_{i \mu}} \hsx \cap \hsx \delta_{i \mu} N_{i \mu}.
\]
The cardinality of the $\Gamma_\gamma$-conjugacy class of \mP in $\sP(L_0(\gamma))$ is 
\[
[\Gamma_\gamma : \Gamma_\gamma \cap P]
\]
%%----------------------------------------------------------------------------------------------05
or still
\[
[n_{i \mu} \Gamma_{\gamma_{i \mu}} n_{i \mu}^{-1} : 
n_{i \mu} \Gamma_{\gamma_{i \mu}} n_{i \mu}^{-1} \cap M_{i \mu}].
\]
Since
\[
\frac
{
[n_{i \mu} \Gamma_{\gamma_{i \mu}} n_{i \mu}^{-1} : 
n_{i \mu} \Gamma_{\gamma_{i \mu}} n_{i \mu}^{-1} \cap M_{i \mu}]
}
{
[n_{i \mu} \Gamma_{\gamma_{i \mu}} n_{i \mu}^{-1} : 
n_{i \mu} \Gamma_{\gamma_{i \mu}} n_{i \mu}^{-1} \cap M_0(\delta_{i \mu})]
}
\]
\[
= \ \frac
{1}
{
[n_{i \mu} \Gamma_{\gamma_{i \mu}} n_{i \mu}^{-1} \cap M_{i \mu} : 
n_{i \mu} \Gamma_{\gamma_{i \mu}} n_{i \mu}^{-1} \cap M_0(\delta_{i \mu})]
},
\]
$C(P_{i \mu},\{\delta_{i \mu}\}_{\Gamma_{M_{i \mu}}})$ 
is equal to 
\[
\sum \frac
{
[(\Gamma_{M_{i \mu}})_{\delta_{i \mu}} \cap M_0(\delta_{i \mu}) : 
n_{i \mu} \Gamma_{\gamma_{i \mu}} n_{i \mu}^{-1} \cap M_0(\delta_{i \mu})]
}
{
[n_{i \mu} \Gamma_{\gamma_{i \mu}} n_{i \mu}^{-1} \cap M_{i \mu} : 
n_{i \mu} \Gamma_{\gamma_{i \mu}} n_{i \mu}^{-1} \cap M_0(\delta_{i \mu})]
}
\]
\[
\times \ \frac{1}{\abs{\det(\Ad(\gamma_{i \mu}) | \mfn_{i \mu} -1}}
\]
or still 
\[
\frac{1}
{[(\Gamma_{M_{i \mu}})_{\delta_{i \mu}} : (\Gamma_{M_{i \mu}})_{\delta_{i \mu}} \cap M_0(\delta_{i \mu})]}
\]
\[
\times \ \sum 
[(\Gamma_{M_{i \mu}})_{\delta_{i \mu}} : n_{i \mu} \Gamma_{\gamma_{i \mu}} n_{i \mu}^{-1} \cap M_{i \mu}]
\]
\[
\times \ \frac{1}{\abs{\det(\Ad(\delta_{i \mu})| \mfn_{i \mu} -1}},
\]
the summation in either case extending over all nonempty
\[
\Delta_0(P)_{P_{i \mu}} \cap \delta_{i \mu} N_{i \mu}.
\]
It remains only to show that the second sum is one.  For this, it will be enough to check that a given nonempty
\[
\Delta_0(P)_{P_{i \mu}} \cap \delta_{i \mu} N_{i \mu}
\]
includes
\[
[(\Gamma_{M_{i \mu}})_{\delta_{i \mu}} : n_{i \mu} \Gamma_{\gamma_{i \mu}} n_{i \mu}^{-1} \cap M_{i \mu}]
\]
$\Gamma \cap N_{i \mu}$-conjugacy classes, there being precisely
\[
\abs{\det(\Ad(\delta_{i \mu})| \mfn_{i \mu} -1}
\]
of these.  
Since the problem is one of counting and therefore, in an obvious sense, is 
``conjugation invariant'', we can assume without loss of generality that $\gamma_{i \mu} = \delta_{i \mu}$, 
then drop the $i$ and the $\mu$ from the notation and finish by proving:

\vspace{0.3cm}

%%----------------------------------------------------------------------------------------------06
{\small\bf Lemma 6.2} \ %02
The $\Gamma \cap P$-conjugacy class of $\gamma$ intersected with $\gamma N \cap \Gamma$ includes
\[
[(\Gamma_{M(\gamma)})_\gamma : \Gamma_\gamma \cap M(\gamma)]
\]
$\Gamma \cap N$-conjugacy classes in $\gamma N \cap \Gamma$.

\vspace{0.1cm}

{\small\bf Proof.} \ 
From the definitions, 
\[
\{\gamma\}_{\Gamma \cap P} \hsx \cap \hsx (\gamma N \cap \Gamma) 
\hsx = \hsx 
\{\eta \gamma \eta^{-1} : \eta \in (\Gamma_{M(\gamma)})_\gamma \bullet N \cap \Gamma\}.
\]
And
\[
(\Gamma_{M(\gamma)})_\gamma \bullet N \cap \Gamma \supset
(\Gamma_\gamma \cap P) \bullet (\Gamma \cap N).
\]
Conjugation by the latter gives a $\Gamma \cap N$-conjugacy class.  Suppose that
\[
\eta_1, \eta_2 \in (\Gamma_{M(\gamma)})_\gamma \bullet N \cap \Gamma
\]
and suppose that $\eta_1 \gamma \eta_1^{-1}$ and $\eta_2 \gamma \eta_2^{-1}$ are $\Gamma \cap N$-conjugate, say
\[
\eta_1 \gamma \eta_1^{-1} \hsx = \hsx n_2 \eta_2 \gamma \eta_2^{-1} n_2^{-1}.
\]
Then
\begin{align*}
\gamma \ 
&=\ \eta_1^{-1} n_2 \eta_2 \gamma \eta_2^{-1} n_2^{-1} \eta_1\\
&\implies \eta_1^{-1} n_2 \eta_2 \in \Gamma_\gamma\\
&\implies \eta_1(\Gamma_\gamma \cap P) = n_2 \eta_2 (\Gamma_\gamma \cap P)\\
&\implies \eta_1(\Gamma_\gamma \cap P) \bullet (\Gamma \cap N)
= \eta_2 (\Gamma_\gamma \cap P) \bullet (\Gamma \cap N).
\end{align*}
Because each such coset fills out a full $\Gamma \cap N$-conjugacy class, the number of $\Gamma \cap N$-conjugacy classes 
in question is 
\[
[(\Gamma_{M(\gamma)})_\gamma \bullet N \cap \Gamma : 
(\Gamma_\gamma \cap P)  \bullet (\Gamma \cap N)].
\]
The map 

\[
\begin{cases}
\ \phi:\Gamma \cap P \ra \Gamma_{M(\gamma)}\\
\ \phi(\eta) = \delta, \quad \{\delta\} = M \cap \eta N
\end{cases}
\]
has kernel $\Gamma \cap N$, so our number is 
\[
[\phi((\Gamma_{M(\gamma)})_\gamma \bullet N \cap \Gamma) : 
\phi((\Gamma_\gamma \cap P)  \bullet (\Gamma \cap N))]
\]
or still
\[
[(\Gamma_{M(\gamma)})_\gamma : \Gamma_\gamma \cap M(\gamma)],
\]
%%----------------------------------------------------------------------------------------------07
as desired.

\vspace{0.3cm}

Hence:

\vspace{0.3cm}

{\small\bf Theorem 6.3} \ %03
$\bK(\bH: \alpha : \Gamma_{\sC_0,\sC_0})$ is equal to 
\[
\sum\limits_{i, \mu}
\ 
\sum\limits_{\{\delta_{i \mu}\}_{\Gamma_{M_{i \mu}}} : 
\delta_{i \mu} \in (\Gamma_{M_{i \mu}})_{\sC_0^\dagger,\sC_0^\dagger}}
\frac{1}
{[(\Gamma_{M_{i \mu}})_{\delta_{i \mu}} : (\Gamma_{M_{i \mu}})_{\delta_{i \mu}} \cap M_0(\delta_{i \mu})]}
\]
\[
\times p(\Gamma : P_{i \mu} : \bH(P_{i \mu}) - \bH_\bO(P_{i \mu})) 
\]
\[
\times \int_{M_{i \mu} / A_0^\dagger(\delta_{i \mu}) 
\bullet 
((\Gamma_{M_{i \mu}})_{\delta_{i \mu}} \cap M_0(\delta_{i \mu}))}
\alpha^{P_{i \mu}}(m\delta_{i \mu} m^{-1}) v_{\mfa_0^\dagger(\delta_{i \mu})} 
(m : I_{M_{i \mu}} (\bH_\bO)) d_{M_{i \mu}} (m).
\]

Denote by 
\[
\bCon - \bCl(\bH: \alpha : \Gamma_{\sC_0,\sC_0})
\]
the result of setting $\bH_\bO = \bH$ in Theorem 6.3.  Since
\[
p(\Gamma : P_{i \mu} :0) \hsx = \hsx 0
\]
unless $P_{i \mu} = G$, 
\[
\bCon - \bCl(\bH: \alpha : \Gamma_{\sC_0,\sC_0})
\]
can also be explicated by Theorem 5.4.  It is then possible to go further but the discussion is combinatorially messy and 
unenlightening, thus will be omitted.

\vspace{0.1cm}
%%%%%%%%%%%%%%%%%%%%%%%%%%%%%%%%%%%%%%
%%%%%%%%%%%%%%%%%%%%%%%%%%%%%%%%%%%%%%
%%%%%%%%%%%%%%%%%%%%%%%%%%%%%%%%%%%%%%

\setcounter{page}{1}
\renewcommand{\thepage}{References-\arabic{page}}
\begingroup
\center {\textbf{REFERENCES}}\\
\endgroup
\vspace{0.5cm}
%\[
%\textbf{References}
%\]

\noindent Arthur, J.: 

[1-(a)] \quad
{The characters of discrete series as orbital integrals}, 
\emph{Inv. Math.}
\textbf{32} (1976), 
205-261.

[1-(b)] \quad
{A trace formula for reductive groups I}, 
\emph{Duke Math. J}
\textbf{45} (1978), 
911-952.\\

\vspace{0.1cm}

\begingroup
\noindent {Osborne, M. S., and Warner, G.:}
\endgroup

[2-(a)] \quad
{The Selberg trace formula I}, 
\emph{Crelle's J.}
\textbf{324} (1981), 
1-113.
%%%%%%%%%%%%%%%%%%%

[2-(b)] \quad
{The Selberg trace formula II}, 
\emph{Pacific J. Math.}
\textbf{106} (1983), 
307-496.
%%%%%%%%%%%%%%%%%%%

[2-(c)] \quad
{The Selberg trace formula III}, 
\emph{Memoirs Amer. Math. Soc.}
\textbf{283} (1983), 
1-209.
%%%%%%%%%%%%%%%%%%%

[2-(d)] \quad
{The Selberg trace formula IV}, 
\emph{SLN} 
\textbf{1024} (1983), 
112-263.
%%%%%%%%%%%%%%%%%%%

[2-(e)] \quad
{The Selberg trace formula V}, 
\emph{Trans. Amer. Math. Soc.}
\textbf{286} (1984), 
351-376.
%%%%%%%%%%%%%%%%%%%

[2-(f)] \quad
{The Selberg trace formula VI}, 
\emph{Amer. J. Math.}
\textbf{107} (1985), 
1369-1437.
%%%%%%%%%%%%%%%%%%%

[2-(g)] \quad
{The Selberg trace formula VII}, 
\emph{Pacific J. Math.}
\textbf{140} (1989), 
263-352.
%%%%%%%%%%%%%%%%%%%

[2-(h)] \quad
{The Selberg trace formula VIII}, 
\emph{Trans. Amer. Math. Soc.}
\textbf{324} (1991), 
623-653.

%\begin{thebibliography}{X}

%\bibitem[Arthur-(a)]{Art}
%Arthur, J., 
%\emph{The characters of discrete series as orbital integrals}, 
%Inv. Math. 
%\textbf{32}, (1976), 
%205-261.

%\end{thebibliography}

%\newpage
%\setcounter{page}{1}
%\renewcommand{\thepage}{Index-\arabic{page}}
%\printindex
\end{document}